\documentclass[11pt,oneside,a4paper,reqno]{amsart}
\usepackage[utf8]{inputenc}
\RequirePackage{fix-cm} 
\usepackage[T1]{fontenc}
\usepackage[english]{babel}

\usepackage{lmodern}
\usepackage{amsmath}
\usepackage{amsthm}
\usepackage{amssymb}
\usepackage{mathrsfs}
\usepackage{graphicx}
\usepackage{color}
\usepackage{longtable}
\usepackage{mathtools}
\usepackage{leftidx}
\usepackage{eurosym}
\usepackage{enumerate}
\usepackage{dsfont}
\usepackage{float}
\usepackage[dvipsnames]{xcolor}
\usepackage{todonotes}
\usepackage{enumitem}

\usepackage{hyperref}
\hypersetup{
colorlinks,
citecolor=Mahogany,
linkcolor=Mahogany,
urlcolor=Mahogany}

\makeatletter
\def\@seccntformat#1{\hspace*{0mm}%
 \protect\textup{\protect\@secnumfont
   \ifnum\pdfstrcmp{subsection}{#1}=0 \bfseries\fi
   \csname the#1\endcsname
   \protect\@secnumpunct
     }%
}
 
\usepackage{geometry}
\newcommand{\N}{\mathbb{N}}
\newcommand{\Z}{\mathbb{Z}}

\newcommand{\R}{\mathbb{R}}
\newcommand{\Sp}{\mathbb{S}}
\newcommand{\dl}{\,\mathrm{d}}

\newcommand{\sph}{\mathbb{S}}
\newcommand{\p}{\partial}
\newcommand{\sequence}[2]{(#1_{#2})_{#2}}

\newcommand{\norm}[2]{\|#1\|_{#2}}
\newcommand{\sfT}{\mathsf{T}}

\newcommand{\veps}{\varepsilon}
\newcommand{\Om}{\Omega}

\newcommand{\ol}[1]{\overline{#1}}

\renewcommand{\ae}{a.\,e.\ }

\DeclareMathOperator{\loc}{loc}

\DeclareMathOperator{\dist}{dist}

\DeclareMathOperator{\id}{id}
\DeclareMathOperator{\Tr}{Tr}
\DeclareMathOperator{\conv}{conv}
\DeclareMathOperator{\SL}{SL}
\DeclareMathOperator{\SO}{SO}

\DeclareMathOperator{\VMO}{VMO}

\newtheorem{Theorem}{Theorem}[section]

\newtheorem{Corollary}[Theorem]{Corollary}
\newtheorem{Proposition}[Theorem]{Proposition}
\newtheorem{Lemma}[Theorem]{Lemma}
\newtheorem*{Conjecture}{Conjecture}

\theoremstyle{definition}
\newtheorem{Definition}[Theorem]{Definition}
\newtheorem{Remark}[Theorem]{Remark}
\newtheorem{Setting}[Theorem]{Convention}
\newtheorem{Assumption}[Theorem]{Assumption}

\newenvironment{theorem}[1][]
{ \begin{Theorem}[#1]}
	{\end{Theorem}}

\newenvironment{corollary}[1][]
{ \begin{Corollary}[#1]}
	{\end{Corollary} }

\newenvironment{proposition}[1][]
{ \begin{Proposition}[#1]}
	{\end{Proposition}}

\newenvironment{lemma}[1][]
{ \begin{Lemma}[#1]}
	{\end{Lemma}}

\newenvironment{remark}[1][]
{ \begin{Remark}[#1]}
	{\end{Remark}}

\numberwithin{equation}{section}

\title[]{On the existence of extensions\\for manifold-valued Sobolev maps\\on perforated domains}

\author[C. Gavioli]{Chiara Gavioli}
\author[L. Happ]{Leon Happ}
\author[V. Pagliari]{Valerio Pagliari}

\date{}

\AtEndDocument{
	\bigskip
	\footnotesize{
		\noindent
		{\sc C. Gavioli}, Faculty of Civil Engineering, Czech Technical University, Th\'akurova 7, 16629 Prague 6, Czech Republic.\\		
		E-mail: \href{mailto:chiara.gavioli@cvut.cz}{\tt chiara.gavioli@cvut.cz}
		
		\medskip
		
		\noindent
		{\sc L. Happ}, Departamento de Matem\'aticas, 
            Universidad Aut\'onoma de Madrid, 
            C/ Francisco Tomás y Valiente 7,
            28049 Madrid, Spain,
            and Institute of Analysis and Scientific Computing,
		TU Wien, Wiedner Hauptstra\ss e 8-10, 1040 Vienna, Austria}\\
            E-mail: \href{mailto:leon.happ@uam.es}{\tt leon.happ@uam.es}
            
           \medskip
		
		\noindent
		{\sc V. Pagliari}, Institute of Analysis and Scientific Computing, TU Wien, Wiedner Hauptstra\ss e 8-10, 1040 Vienna, Austria.\\
		E-mail: \href{mailto:valerio.pagliari@tuwien.ac.at}{\tt valerio.pagliari@tuwien.ac.at}
}

\begin{document}
	
	\begin{abstract}
		Motivated by manifold-constrained homogenization problems, we construct suitable extensions for Sobolev functions defined on a perforated domain and taking values in a compact, connected $C^2$-manifold without boundary. The proof combines a by now classical extension result for the unconstrained case with a retraction argument that hea\-vi\-ly relies on the topological properties of the manifold.
        With the ultimate goal of providing necessary conditions for the existence of extensions for Sobolev maps between manifolds, we additionally investigate the relationship between this problem and the surjectivity of the trace operator for such functions.
		
		\medskip
		
		\noindent
		{\it 2020 Mathematics Subject Classification:}
        54C20, 
        46T10, 
		46E35, 
        55S35 
		
		\smallskip
		\noindent
		{\it Keywords and phrases:}
		Extension of functions, perforated domain, manifold constraint, trace operator
		
	\end{abstract}
	
	\maketitle
	
	{\parskip=0em \tableofcontents}
	

\section{Introduction}\label{sec:intro}

In homogenization problems on perforated domains, a crucial step is the formulation of a meaningful notion of convergence that accounts for the fact that the functions involved are defined on changing domains. A possible approach amounts to extending the functions to the holes with the help of extension operators preserving certain boundedness properties independently of the scale of the microperforations, and to studying the convergence of these extended functions, which are now all defined on the same domain. This procedure has been effectively used in the literature: we refer to the prominent paper by Acerbi et al. \cite{AcPiMaPe92} and the references therein (in particular \cite{CioSJP79}, where for the first time an extension operator for Sobolev functions on perforated domains was constructed), as well as to the monograph \cite{BrDe98}.

In the references mentioned, the functions under consideration can freely range in a certain Euclidean space $\R^l$.
In physics, however, problems that are formulated in terms of manifold-valued Sobolev functions arise naturally, e.g., in micromagnetics with $\Sp^2$ \cite{Br63,HuSch08}, in finite-strain plasticity with $\SL(3)$ \cite{CaHaMie02,MaMi09,DaGaPa1}, or in the modeling of rigid body rotations with $\SO(3)$ \cite{Chris20,Engl22}. To our knowledge, no result analogous to the one in \cite{AcPiMaPe92} is available for such scenarios. With the aim of filling this gap,
the present paper investigates the existence of suitable
extensions of manifold-valued Sobolev maps defined on a perforated domain $\Omega_\veps$, that is, loosely speaking,
the region obtained by periodically removing from an open set $\Omega\subset \R^d$ holes whose size is of the order $\veps>0$.
To be precise, for $\veps>0$, \(1<p<\infty\), and a given compact, connected $C^2$-manifold $N$ without boundary, we are looking for a function $\tilde{f}_\veps \in W^{1,p}(\Om;N)$ such that $\tilde{f}_\veps|_{\Om_\veps} = f$, and that additionally admits
bounds of the form
\begin{gather}
	\norm{\tilde{f}_\veps}{L^p(\Om;N)} \le C\norm{f}{L^p(\Om_\veps;N)},\label{eq:intro_bounds_lp} \\
	\norm{D\tilde{f}_\veps}{L^p(\Om;\R^{l\times d})} \le C\norm{Df}{L^p(\Om_\veps;\R^{l\times d})},\label{eq:intro_bounds_grad_lp}
\end{gather}
for some constant $C>0$ independent of the scale $\veps$.

Provided $\Om_\veps$ has a certain structure,
our main result in this direction is Theorem~\ref{thm:target_ext}.
There, we exhibit an extension fulfilling \eqref{eq:intro_bounds_lp}--\eqref{eq:intro_bounds_grad_lp} under the assumptions that \(p<d\) and that, denoting by $\lfloor \,\cdot\, \rfloor$ the floor function, the first \(\lfloor p-1\rfloor\) homotopy groups of \(N\) are trivial.
The latter is a topological condition called {\em \(\lfloor p-1\rfloor\)-connectedness} (see Section~\ref{subsec:homotopy} for the definitions).
Such hypotheses are dictated by our approach to the proof, which is inspired by that of \cite[Theorem~6.2]{HaLi87}, where the existence of finite-energy extensions of Sobolev mappings between manifolds is proved in the subcritical case $p<d$.
To be precise, our strategy combines the extension result for the unconstrained case in~\cite{AcPiMaPe92}
(see Theorem~\ref{thm:reflect_ext} and Corollary~\ref{cor:reflect})
with a retraction argument from \(\R^l\) onto the manifold \(N\) proposed in~\cite{HaLi87} (see Lemma~\ref{lem:lip_retract}).
Indeed, by appealing to Theorem~\ref{thm:reflect_ext} or Corollary~\ref{cor:reflect},
an extension can be found that takes values in the whole space $\R^l$.
Then, the conclusion follows by a careful composition with a suitable retraction on the target manifold $N$.
The \(\lfloor p-1\rfloor\)-connectedness condition ensures that, even though in general the retraction will possess singularities in the form of lower dimensional planes, it will still enjoy \(W^{1,p}\)-regularity.

Summarizing, when $p<d$, a sufficient condition for the existence of extensions for Sobolev functions taking values in a sufficiently smooth compact manifold $N$ without boundary is that \(N\) is \(\lfloor p-1\rfloor\)-connected.
On the other hand, we cannot provide a satisfactory answer about the existence of an extension if $N$ is not \(\lfloor p-1\rfloor\)-connected, and in general, the problem of devising necessary conditions for the existence of extensions as in \eqref{eq:intro_bounds_lp}--\eqref{eq:intro_bounds_grad_lp} is at the moment
mostly open.
Here, we contribute to its understanding in Section~\ref{sec:conn_trace} by analyzing the extension problem in the simpler case of a single constituting cell,
and by thoroughly checking the conditions that arise from the study of the surjectivity of the trace operator for Sobolev functions between two compact Riemannian manifolds, both for the cases $p\ge d$ (see Proposition~\ref{prop:gen_nec_cond}) and $p<d$ (see Propositions~\ref{prop:gen_nec_cond_smallp} and \ref{prop:case2.2}).

Several works have already been devoted to the problem of finding sufficient and necessary conditions for the surjectivity of the trace operator for Sobolev maps between a manifold with boundary and a given target manifold $N$ (cf.\ \cite{HaLi87, BeDe95, Be14, MiVS21, MaVS23, VS24}).
In the case $p<d$, the assumption that $N$ is $\lfloor p-1 \rfloor$-connected plays a pivotal role: it is now well-known that $\lfloor p-1 \rfloor$-connectedness guarantees that the trace operator is surjective (see \cite{HaLi87}). On the other hand, the recent preprint \cite{VS24} suggests that, when \(p<d\), a sufficient and necessary condition for the surjectivity of the trace operator 
that is weaker than \(\lfloor p-1\rfloor\)-connectedness can be devised.
When \(p>d\), it follows from the Sobolev embeddings that the surjectivity of the Sobolev trace operator is equivalent to the same property in the space of continuous functions. It turns out that the same holds true for the critical case \(p=d\) (see \cite[Theorems~1 and 2]{BeDe95}).
The extension problem for continuous functions has been extensively studied in the context of algebraic topology, and for further reading on this point we refer the reader to the monographs \cite{Hu59,Ha02}.

In addition to its intrinsic value, the importance of the present extension result lies in its potential applications in homogenization theory, on which there is nowadays a large body of work (see, e.g., \cite{BrDe98,CioDo99,BeRy18}) including the case of manifold-valued Sobolev maps (see \cite{BaMi10}).
In particular, when dealing with perforated domains, an extension that preserves a given target manifold condition can be used to perform the asymptotic analysis
of specific physical models, as we briefly show in Section~\ref{sec:ex} in the context of micromagnetics.
Lastly, we mention that the extension problem also plays an important role in the study of high-contrast homogenization models, namely, when the holes in the perforated domain are filled with a material that has drastically different properties than the matrix substance
(see \cite{CC12,CCNe17,DaKrPa,DaGaPa2}).
The possibility of using the extension result given by Theorem~\ref{thm:target_ext} in a high-contrast homogenization problem in micromagnetics was investigated in \cite{DaHa24}.


\subsection*{Outline of the paper}

In Section~\ref{sec:prelim} we fix notations and definitions, and recall some notions from geometry and algebraic topology that we use throughout the paper. We also briefly compare different definitions of perforated domains that ultimately influence our main result on the existence of extensions (cf.\ Theorem~\ref{thm:target_ext} and Corollary~\ref{cor:target_ext}).
Thereafter, Section~\ref{sec:constr_ext} is devoted to the proof of Theorem~\ref{thm:target_ext}.
The two main ingredients,
the extension operator from \cite{AcPiMaPe92} and the retraction map from~\cite{HaLi87}, are recalled.
We additionally show how such an extension can be advantageously employed in a model from micromagnetics.
Shifting the focus to necessary conditions for the existence of extensions,
in Section~\ref{sec:conn_trace} we collect some results on the surjectivity of the trace operator for Sobolev maps between manifolds and establish connections to our central question.
Finally, in Section~\ref{sec:concl} we draw conclusions and present possible further applications.


\section{Preliminaries}\label{sec:prelim}


\subsection{Notation and mathematical setting}

Let $d,l\in \N$, $d,l\geq2$.
The symbol $|\,\cdot\,|$ is adopted without further specifications for the Euclidean norms in $\R^d$, $\R^l$, and $\R^{l\times d}$,
and the symbol $\dist$ for the associated distances.
If $A \subset \R^d$ is a measurable set, we denote by $\mathcal{L}^d(A)$ its $d$-dimensional Lebesgue measure,
while with \(\mathcal{H}^{d-1}\) we refer to the \(d-1\)-dimensional Hausdorff measure.

For \(\lambda>0\) and a set \(A\subset \R^d\), we define \(A(\lambda)\coloneqq \{x\in A:\ \dist(x,\p A)>\lambda\}\), that is,
$A(\lambda)$ is a retracted set of $A$.
We set \(Q\coloneqq \left(-\frac{1}{2},\frac{1}{2} \right)^d\subset \R^d\) for the unit cube with center at the origin.
If not specified differently, the letter \(C\) stands for a positive constant that may change from line to line. 

Let $N$ be a compact $C^2$-submanifold of $\R^l$ without boundary.
In Section~\ref{sec:conn_trace}, we also consider the Riemannian structure on $N$ inherited from the ambient space.
Note that $N$, for a sufficiently large $l$, can also be the realization
of an abstract compact $C^2$ manifold $\tilde{N}$
obtained by means of the Whitney embedding theorem (cf.\ \cite{Wh44}), or by means of the Nash embedding theorem in the Riemannian case (cf.\ \cite{Na56}).

For \(\delta>0\), we define the \emph{tubular neighborhood} \(N_\delta\) of \(N\) as
\begin{equation*}
	N_\delta
	\coloneqq\{x\in \R^l:\, \dist(x,N)< \delta\}.
\end{equation*}
Under our regularity assumptions on \(N\), there exists some particular \(\delta>0\) such that
for each \(x\in N_\delta\) there is a unique point \(\xi(x)\in N\)
with the property that \(\dist(x,N) = |x-\xi(x)|\).
The map \(\xi:N_\delta\to N\) is called the {\it nearest point projection} and loses one degree of regularity with respect to \(N\)
(cf.\ \cite[Theorem~6.17 and Proposition~6.18]{Lee03} as well as \cite{Fo84}),
so that, in our setting, the next statement holds true:

\begin{lemma}\label{lem:nearest_pp}
    Let $N$ be a compact $C^2$-submanifold of $\R^l$ without boundary.
	There exists \(\delta>0\) such that
    the nearest point projection \(\xi:N_\delta\to N\) is well-defined and
    fulfills \(\xi\in C^1(\overline{N_\delta};N)\).
\end{lemma}
In what follows, the symbol $\delta$ refers always to the parameter in the previous lemma.

Finally, the functional setting of our analysis is the the space $W^{1,p}(\Om;N)$
of $N$-valued Sobolev functions on $\Om$.
Specifically, for $1<p < \infty$, this space is defined as follows:
\[
	W^{1,p}(\Om;N) \coloneqq \{ f \in 	W^{1,p}(\Om;\R^l) :\ f(x)\in N \text{ for \ae } x\in \Om\}.
\]


\subsection{Reminders on homotopy groups and a retraction result}\label{subsec:homotopy}

We recall here some notions from homotopy theory.
We refer to \cite{Ha02,Hu59} for a thorough treatment of the topic.
	
Let $X$ be a path-connected topological space and let us arbitrarily fix $x_0\in X$.
We denote by $Q^i$ the $i$-dimensional unit cube $\left(-\frac{1}{2},\frac{1}{2} \right)^i$.
For $i\in \N\setminus\{0\}$, the {\it $i$-th homotopy group of $X$}, $\pi_i(X,x_0)$, is the set of all homotopy classes of the continuous maps $f\colon \ol{Q}^i \to X$ such that $f(\partial Q^i)=\{x_0\}$, where $\partial Q^i$ is the topological boundary of $Q^i$.
Recall that under the assumption of path-connectedness of $X$ any two groups \(\pi_i(X,x_0)\) and \(\pi_i(X,x_1)\) are isomorphic, and thus the fundamental group does not depend on the choice of $x_0$ (cf.\ \cite[Proposition~1.5 and pp.~341]{Ha02}).
Hence, to unburden the notation, we just write $\pi_i(X)$ if $X$ is path-connected.

The space \(X\) is said to be \emph{\(n\)-connected} if \(\pi_i(X)=0\) for all \(i\le n\) (cf.\ \cite[pp.~346]{Ha02}), that is,
if for all \(i\le n\) the \(i\)-th homotopy group consists only of the trivial homotopy class.
Denoting by $\sph^i$ the $i$-sphere of radius $1$, it follows from the definition of homotopy groups that this is equivalent to the property that every continuous map \(\sph^i\to X\) can be contracted to a point (i.e., it is homotopic to a constant map).
As a consequence of the fact that the pair \((B^{i+1},\sph^i)\) has the {\it homotopy extension property} (cf.\ \cite[Chapter~I, Proposition~4.2]{Hu59} and \cite[Proposition~0.16]{Ha02}),
it also turns out that 
\(n\)-connectedness is equivalent to the property that for all \(i\le n\)
every continuous map \(\sph^i\to X\) can be extended to a continuous map \(B^{i+1}\to X\),
where $B^{i+1}$ denotes the \((i+1)\)-dimensional unit ball.

The features we have just recalled underpin the following result,
which provides a retraction from the whole space to the target manifold and
will come in handy in the proof of Theorem~\ref{thm:target_ext}.
The statement is essentially the one in \cite[Lemma~2.2]{BPV14},
which is in turn a reformulation of \cite[Lemma~6.1]{HaLi87}.

\begin{lemma}\label{lem:lip_retract}
    Let $N$ be a compact $C^2$-submanifold of \(\R^l\),
    and let $\delta >0$ be the parameter provided by Lemma \ref{lem:nearest_pp}.
    Let also $R>0$ be such that
    the $\delta$-tubular neighborhood \(N_\delta\) of $N$
    is contained
    in the $l$-dimensional cube $Q_R \coloneqq (-R,R)^l$,
	and let finally $j \in \{0, \dots, l-2\}$. If 
	\begin{equation}\label{eq:pi}
		\pi_0(N) = \dots =\pi_j(N) = 0,
	\end{equation}
	then there exist a closed subset $X \subset Q_R \setminus N$ contained in a finite union of (\(l-j-2\))-dimensional planes and a locally Lipschitz map \(P : Q_R\setminus X \to N\) such that the following hold:
	\begin{enumerate}[label=(\roman*)]
		\item \label{item:retraction}
			$P$ is a retraction onto $N$, that is, $P(y)=y$ for every $y \in N$;
		\item\label{item:inclusion}
			on \(N_\delta\) the map \(P\) is equal to the nearest point projection, and, in particular, one has
			\begin{equation*}
				N_\delta\cap X
				= \emptyset;
			\end{equation*}
		\item\label{item:projection_grad_lp}
				there exists a constant $C>0$ depending on $l$, $N$, and \(R\) such that, for every $y \in Q_R\setminus X$,
			$$
				|DP(y)| \le \frac{C}{\dist(y,X)},
			$$
				and, in particular, \(DP \in L^p(Q_R;\R^{l\times d})\) for $1\le p<j+2$.
	\end{enumerate} 
\end{lemma}
In Appendix~\ref{sec:retract_lem} we give an account of the construction of the map $P$, and we make explicit the role of condition \eqref{eq:pi}.


\subsection{Perforated domains and available extension results}\label{sec:recap_extensions}

Throughout the paper $\Omega \subset \R^d$ denotes a bounded open set.
As we mentioned in the introduction,
in the present contribution we deal with functions defined
on the so called \emph{perforated domains}, i.e.,
on sets obtained from $\Om$ by removing small holes from it.
In this subsection, we first revise two possible definitions of perforated domain,
then we recall the corresponding known extension results for (unconstrained) Sobolev maps
(see Theorem~\ref{thm:reflect_ext} and Corollary~\ref{cor:reflect}).

The starting point to define a perforated domain is
an open, connected set \(E\subset \R^d\)
that has Lipschitz boundary and is periodic, i.e.,
it fulfills $E+e_i = E$ for all $i=1,\dots, d$,
where \((e_i)_{i=1,\dots d}\) is the canonical basis of \(\R^d\).
We use the set $E$ to describe the microstructure of the perforated set $\Om_\veps$ as follows.

First, we define
\[
    Q_1\coloneqq Q\cap E
    \quad\text{and}\quad
    Q_0\coloneqq Q\setminus \overline{Q_1}.
\]
The sets $Q_0$ and $Q_1$ represent, respectively,
the perforation and the ``solid'' part of the unit cell $Q$ at the scale $\veps = 1$.
For our subsequent analysis (see the proof of Lemma \ref{lemma:estimate_measure}),
we also need to require that
\begin{equation}\label{eq:meas_q1}
   0 < \mathcal{L}^d(Q_1) < 1. 
\end{equation}

Second, we translate and rescale $Q_0$ to define the set of micropores and the perforated set.
Precisely, recalling that,
for \(\lambda>0\) and \(A\subset \R^d\), $A(\lambda)$ is a retracted set of $A$,
we introduce the family of translations
\begin{equation}\label{eq:second_setup1}
    Z_\veps\coloneqq \{z\in \Z^d : \veps (Q_0+z) \subset \Om(\veps\lambda) \},
\end{equation} 
where $\lambda>0$ is a given, fixed parameter.
Then, as in~\cite{DaGaPa2} and similarly to, e.g., \cite[Chapter~4]{OSY12},
we define the perforations \(\Om_{0,\veps}\) and the perforated domain \(\Om_\veps\) of scale \(\veps>0\) respectively as 
\begin{equation}\label{eq:second_setup2}
\Om_{0,\veps}
\coloneqq \bigcup_{z\in Z_\veps} \veps(Q_0+z)
\quad \text{ and }\quad
\Om_\veps \coloneqq \Om\setminus \overline{\Om_{0,\veps}},
\end{equation}
with the understanding that \(\Om_{0,\veps}=\emptyset\) for \(Z_\veps=\emptyset\).
Note that, in particular, it holds \(\mathcal{L}^d(\Om_\veps)>0\)
whenever \(Z_\veps\neq\emptyset\),
because \( \Om_\veps \supset \Om \setminus \Om(\veps\lambda)\).

Some comments on the positions in \eqref{eq:second_setup1}--\eqref{eq:second_setup2} are in order.
The subfamily $Z_\veps$ of integer-valued translations grants that
the corresponding holes are well-separated from the boundary of \(\Om\),
or, in other words, that
the perforations \(\Om_{0,\veps}\) are compactly contained in \(\Om\),
even though -- in contrast to the seminal works \cite{CioSJP79,CoDo88,DoMo90} -- we do not assume that
the unit perforation $Q_0$ is compactly contained in $Q$ ($Q_0\Subset Q$).
From this perspective,
the notion of perforated domain in \eqref{eq:second_setup1}--\eqref{eq:second_setup2}
stands somehow in between the restrictive scenario when $Q_0\Subset Q$ and
the more general definition from~\cite{AcPiMaPe92},
which was further employed in \cite{BraChi94,BraChiPia04,BraChiDE21}
and which we recall next. 

For $\veps>0$ and
for the same open, connected, and periodic set \(E\subset \R^d\)
with Lipschitz boundary as above, let us define
\begin{equation}\label{eq:stiff_gen}
	\Om_\veps^\star
	\coloneqq \Om\cap \veps E.
\end{equation}
Note that \(\Om_\veps^\star\) can be also written as
\begin{equation}\label{eq:stiff}
\Om_\veps^\star = \Om \cap \left[ \bigcup_{z\in Z_\veps^\star} \veps(Q_1+z) \right],
\quad
\text{with } Z_\veps^\star \coloneqq \{z\in \Z^d : \veps (Q+z)\cap \Om\ne \emptyset\}.
\end{equation}
The set of microperforations is then obtained as
\begin{equation}\label{eq:soft}
    \Om_{0,\veps}^\star
    \coloneqq \Om\setminus \overline{\Om_\veps^\star}
    = \Om \cap \left[ \bigcup_{z\in Z_\veps^\star} \veps(Q_0+z) \right].
\end{equation}

Note that in both definitions the corresponding microperforations \(\Om_{0,\veps}\) and \(\Om_{0,\veps}^\star\) might be connected.
Instead, a relevant difference between the two scenarios is that $\Om_{0,\veps}^\star$ is allowed to intersect the boundary \(\p \Om\) of $\Om$,
while that was not the case with \(\Om_{0,\veps}\).
Such a circumstance prevents the bounds on the extension operators obtained in~\cite{AcPiMaPe92} from holding up to \(\p \Om\).
Precisely, in \cite[Theorem~2.1]{AcPiMaPe92} the following is shown.

\begin{Theorem}\label{thm:reflect_ext}
	Let \(1\le p<\infty\).
    Let also \(\Om\subset \R^d\) be a bounded open set,
    and for $\veps>0$ let \(\Om_\veps^\star\) be as in \eqref{eq:stiff_gen}.
	Then, there exist a linear and continuous extension operator
	\begin{equation*}
		\sfT_\veps^\star : 
		W^{1,p}(\Om_\veps^\star;\R^l) \to W_{\loc}^{1,p}(\Om;\R^l)
	\end{equation*}
	and some strictly positive constants $\mu=\mu(d,p,E)$, $C=C(d,p,E)$ independent of \(\veps\) and \(\Om\) such that
	\begin{gather}
		\sfT_\veps^\star f = f \quad \text{\ae in } \Om_\veps^\star,\label{eq:ex_op_prop1}\\
		\norm{\sfT_\veps^\star f}{L^p(\Om(\veps \mu);\R^l)} 
			\le C\norm{f}{L^p(\Om_\veps^\star;\R^l)},\label{eq:ex_op_prop2}\\
		\norm{D(\sfT_\veps^\star f)}{L^p(\Om(\veps \mu);\R^{l\times d})} \le C\norm{Df}{L^p(\Om_\veps^\star;\R^{l\times d})},\label{eq:ex_op_prop3}
	\end{gather}
	for every \(f \in W^{1,p}(\Om_\veps^\star;\R^l)\). 
\end{Theorem}

If $\Om$ has Lipschitz boundary and, instead of~\eqref{eq:stiff_gen},
the definition of perforated domain in \eqref{eq:second_setup2} is chosen, that is,
if we restrict to \(Z_\veps\subset Z_\veps^\star\) to describe the microperforations,
a variation of Theorem~\ref{thm:reflect_ext} can be proved
where the estimates on the extension operators hold up to the boundary.

\begin{corollary}\label{cor:reflect}
Let \(1\le p<\infty\). Let also \(\Om\subset \R^d\) be an open bounded set with Lipschitz boundary,
and for $\veps>0$ let the perforated domain \(\Om_\veps\) be as in \eqref{eq:second_setup2}.
Then, there exist a linear and continuous extension operator
\begin{equation*}
		\sfT_\veps: 
		W^{1,p}(\Om_\veps;\R^l) \to W^{1,p}(\Om;\R^l)
\end{equation*}
and a constant \(C=C(d,p,E)>0\) independent of \(\veps\) and \(\Om\) such that 
\begin{gather}
		\sfT_\veps f = f \quad \text{\ae in } \Om_\veps,\label{eq:ex_op_simplified_prop1}\\
		\norm{\sfT_\veps f}{L^p(\Om;\R^l)} \le C\norm{f}{L^p(\Om_\veps;\R^l)},\label{eq:ex_op_simplified_prop2}\\
		\norm{D(\sfT_\veps f)}{L^p(\Om;\R^{l\times d})} \le C\norm{Df}{L^p(\Om_\veps;\R^{l\times d})},\label{eq:ex_op_simplified_prop3}
\end{gather}
for every \(f \in W^{1,p}(\Om_\veps;\R^l)\).
\end{corollary}

Corollary~\ref{cor:reflect} can be established in a fashion analogous to Theorem~\ref{thm:reflect_ext},
the only difference between the two being the geometry of the respective perforated sets, and
therefore we do not include the detailed proof here.
However, for later use in this contribution, we sketch the argument in the following remark.

\begin{remark}\label{rmk:small_diff}
    It is well-known (see e.g.~\cite[Lemma 2.5]{AcPiMaPe92}) that
    a Sobolev function defined on a halfspace may be extended by reflection,
    and that for more general domains with Lipschitz boundary
    the extension is obtained, loosely speaking,
    by locally reflecting the function and then by suitably patching together the local extensions.
    At their core, the proofs of Theorem~\ref{thm:reflect_ext} and of Corollary~\ref{cor:reflect} rely on such standard results,
    and consists in constructing extensions on each of the micropores \(\veps(Q_0+z)\),
    that are then in turn patched together
    by means of a partition of unity subordinated to a covering of the holes.
    As a consequence of this procedure,
    on the inside of the holes and very close to their boundary,
    the extension \(\sfT_\veps f\) or \(\sfT_\veps^\star f\) will only take values in the convex hull of the range of the original function \(f\),
    while away from the boundary it will be equal to the mean value of $f$
    (see the proof of \cite[Lemma~2.6]{AcPiMaPe92}).
    This means that, in the case we are interested in, 
    that is, when we consider functions in \( W^{1,p}(\Om_\veps;N)\) for some compact submanifold \(N\subset \R^l\),
    there holds \(\sfT_\veps f(\Om)\subset \conv(N)\) and \(\sfT_\veps^\star f(\Om)\subset \conv(N)\),
    where \(\conv(N)\) the convex hull of \(N\) in \(\R^l\). 
\end{remark}

We point out that our main result, Theorem~\ref{thm:target_ext}, is tailored to the definition of perforated domain given in \eqref{eq:second_setup1}--\eqref{eq:second_setup2} above.
In this case, by grounding on Corollary~\ref{cor:reflect},
we will construct extensions with \(L^p\)-bounds that hold on the whole set $\Om$.
Then, in Corollary~\ref{cor:target_ext}, we provide an extension result for the more general notion of perforated domain in \eqref{eq:stiff_gen}.
In this setting, the presence of holes that may intersect the boundary combined with the manifold constraint allows us to construct extensions only up to a relatively compact subset of $\Om$,
see the discussion before Corollary~\ref{cor:target_ext}.

We conclude this subsection with a lemma
that will be instrumental in establishing the $L^p$-bound in~\eqref{eq:intro_bounds_lp}.
Indeed, it turns out that our strategy for proving this estimate requires bounding the ratio
between the measure of the perforations and that of the perforated set.
As the following lemma shows, for $\veps>0$ sufficiently small
such a ratio may be controlled uniformly in $\veps$
by a quantity that depends only on the microstructure.

\begin{lemma}\label{lemma:estimate_measure}
    Let $\veps>0$ and let $\Om_{0,\veps}$ and $\Om_\veps$ be as in \eqref{eq:second_setup1}--\eqref{eq:second_setup2}.
    There exists a real number
    \(\ol{\veps} = \ol{\veps}(d,\Om,Q_0)>0\) such that
    \begin{equation}\label{eq:estimate_measure}
       \frac{\mathcal{L}^d(\Om_{0,\veps})}{\mathcal{L}^d(\Om_\veps)}
       \leq \frac{1+\mathcal{L}^d(Q_0)}{\mathcal{L}^d(Q_1)}
       \qquad\text{for every } \veps \le \ol{\veps}.
    \end{equation}
\end{lemma}
\begin{proof}
    Since $\Om_{0,\veps}\subset \Om$ and $\Om_\veps \coloneqq \Om\setminus \overline{\Om_{0,\veps}}$ (see \eqref{eq:second_setup2}),
    we clearly have
    \begin{equation}\label{eq:estimate_measure0}
    \mathcal{L}^d(\Om_\veps) = \mathcal{L}^d(\Om)-\mathcal{L}^d(\Om_{0,\veps}).
    \end{equation}
    In order to estimate the right-hand side,
    in what follows we bound the quantity $\mathcal{L}^d(\Om_{0,\veps})$ from above.
    To this aim, we introduce some further notation.

    First, recalling the definition of $Z_\veps^\star$ in \eqref{eq:stiff},
	we define its disjoint subsets
	\begin{align}
            \begin{split}\label{eq:Zeps}
             Z_\veps^{\star,\textup{int}} &\coloneqq \{z\in Z_\veps^\star : \veps(Q+z) \cap \partial\Om = \emptyset\}, \\
        	    Z_\veps^{\star,\textup{bd}} &\coloneqq
                    \{z\in Z_\veps^\star : \veps(Q+z) \cap \partial\Om \neq \emptyset\},
            \end{split}
        \end{align}
	that is, the subsets of translations in $Z_\veps^\star$ such that
    the corresponding $\veps$-cubes are completely contained in $\Om$ or,
    respectively, such that the corresponding $\veps$-cubes have some intersection with the boundary $\partial\Om$.
    We are interested in the cardinalities of such sets,
    that we denote respectively by
    $ \#Z_\veps^{\star,\textup{int}} $ and $\#Z_\veps^{\star,\textup{bd}} $.
    Second, let us set for shortness
    $\omega \coloneqq \mathcal{L}^d(\Om)$ and $q_0 \coloneqq \mathcal{L}^d(Q_0)$.
	
	Then, by the definition of $\Om_{0,\veps}$ in \eqref{eq:second_setup2} and
    the inclusion $Z_\veps \subset Z_\veps^\star$, we have
    \begin{equation*}
	\mathcal{L}^d(\Om_{0,\veps})
        = q_0 \veps^d \#Z_\veps
        \le q_0 \veps^d \#Z_\veps^\star,
    \end{equation*}
    where $\#Z_\veps$ and $\#Z_\veps^\star$ are the cardinalities of the two sets.
    Since $Z_\veps^\star = Z_\veps^{\star,\textup{int}} \cup Z_\veps^{\star,\textup{bd}}$,
    we further obtain
    \begin{equation}\label{eq:estimate_measure_Om0eps}
	\mathcal{L}^d(\Om_{0,\veps}) \le 
    q_0 \veps^d \Big( \# Z_\veps^{\star,\textup{int}} + \# Z_\veps^{\star,\textup{bd}} \Big).
    \end{equation}

    Our next task is to estimate the cardinalities
    $ \# Z_\veps^{\star,\textup{int}} $ and $\# Z_\veps^{\star,\textup{bd}} $.
    As for the first one, observe that, by definition of $Z_\veps^{\star,\textup{int}}$,
    \[
        \bigcup_{z\in Z_\veps^{\star,\textup{int}}} \veps (Q+z) \subset \Om,
    \]
    whence
    \begin{equation}\label{eq:card_Zint}
        \veps^d \# Z_\veps^{\star,\textup{int}} \leq \omega.
    \end{equation}
    For what concerns $\# Z_\veps^{\star,\textup{bd}}$, we note that
    \[
        \bigcup_{z\in Z_\veps^{\star,\textup{bd}}} \veps (Q+z)
        \subset \{x \in \R^d : \dist(x,\partial\Om) < \veps\sqrt{d}\},
    \]
    and thanks to the regularity of $\partial\Om$, we deduce
    \begin{equation}\label{eq:card_Zbd}
      \veps^d \# Z_\veps^{\star,\textup{bd}}
      \leq 2 \veps \sqrt{d} \mathcal{H}^{d-1}(\partial\Om).
    \end{equation}
    By plugging \eqref{eq:card_Zint} and \eqref{eq:card_Zbd} into \eqref{eq:estimate_measure_Om0eps}, we find
    \begin{equation}\label{eq:estimate_measure_Om0eps_2}
        \mathcal{L}^d(\Om_{0,\veps})
        \leq q_0 \Big( \omega + 2 \veps \sqrt{d} \mathcal{H}^{d-1}(\partial\Om) \Big).
    \end{equation}

    In order to employ \eqref{eq:estimate_measure_Om0eps_2} in \eqref{eq:estimate_measure0}
    to bound $\mathcal{L}^d(\Om_\veps)>0$ from below,
    we impose that
    \[
        \omega - q_0 \Big( \omega + 2 \veps \sqrt{d} \mathcal{H}^{d-1}(\partial\Om) \Big) > 0 ,
    \]
    that is,
    \[
        \veps <
        \frac{\omega}{2 \sqrt{d}\mathcal{H}^{d-1}(\partial\Om)} \frac{1-q_0}{q_0}
    \]
    (observe that the right-hand side is well-defined and strictly positive because of \eqref{eq:meas_q1}).
    Now let, e.g.,
    \begin{equation}\label{eq:threshold}
        \ol\veps 
        \coloneqq \frac{\omega}{4 \sqrt{d}\mathcal{H}^{d-1}(\partial\Om)} \frac{1-q_0}{q_0}.
	\end{equation}
     From \eqref{eq:estimate_measure0} and \eqref{eq:estimate_measure_Om0eps_2} we infer that for every $\veps \le \ol{\veps}$ it holds that 
    \[
    \frac{\mathcal{L}^d(\Om_{0,\veps})}{\mathcal{L}^d(\Om_\veps)}
    \leq \frac{q_0 \Big( \omega + 2 \ol\veps \sqrt{d} \mathcal{H}^{d-1}(\partial\Om) \Big)}{ \omega - q_0 \Big( \omega + 2 \ol\veps \sqrt{d} \mathcal{H}^{d-1}(\partial\Om) \Big)}
    = \frac{1+q_0}{1-q_0},
    \]
    as desired.
\end{proof}

\begin{remark}\label{rmk:suff_dist}
    In view of the proof of Corollary~\ref{cor:target_ext} and in the spirit of Lemma~\ref{lemma:estimate_measure}, we can derive an estimate for the Lebesgue measure of the set $\Om_{0,\veps}^\star\cap\Om(\veps \mu)$, where \(\Om_{0,\veps}^\star\) is as in \eqref{eq:soft} and \(\Om(\veps \mu)\) is the retracted set appearing in Theorem~\ref{thm:reflect_ext}. Let us denote by \(I(\Om(\veps \mu))\) the set of all integer vectors \(z\in \Z^d\) such that \(\veps(Q+z)\cap \Om(\veps \mu) \neq \emptyset\). It follows from \cite[Lemma~2.7]{AcPiMaPe92} that \(\veps(Q+z)\subset \Om\) for all \(z\in I(\Om(\veps \mu))\). Recalling the definition of \(Z_\veps^{\star,\textup{int}}\) in \eqref{eq:Zeps}, we conclude that \(I(\Om(\veps \mu))\subset Z_\veps^{\star,\textup{int}}\). In particular there holds
    \begin{equation}\label{eq:gen_soft_est}
    \mathcal{L}^d\big(\Om_{0,\veps}^\star\cap\Om(\veps \mu)\big)
        \le \sum_{z \in Z_\veps^{\star,\textup{int}}} \veps^d \mathcal{L}^d(Q_0).
    \end{equation}
\end{remark}

	
\section{Construction of an extension when \texorpdfstring{$N$}{N} is \texorpdfstring{$\lfloor p-1\rfloor$}{[p-1]}-connected}\label{sec:constr_ext}

In this section, we prove
the existence of suitable extensions for manifold-valued Sobolev functions defined on a perforated domain,
when the integrability exponent $p$ is smaller
than the dimension $d$
and the manifold $N$ is $\lfloor p-1\rfloor$-connected.
We first address the case of the perforated domain as defined in \eqref{eq:second_setup2}.

\begin{Theorem}\label{thm:target_ext}
    Let \(\Om\subset \R^d\), $d\geq 2$, be an open and bounded
    set with Lipschitz boundary, and
    let the perforated domain \(\Om_\veps\) be defined as in \eqref{eq:second_setup1}--\eqref{eq:second_setup2}.
    Let also $N$ be a compact $C^2$-submanifold of $\R^l$, $l\geq 2$, without boundary.
    
    If \(1< p<d\) and
	\begin{equation}\label{eq:pi_target}
	\pi_0(N) = \dots =\pi_{\lfloor p-1\rfloor}(N) = 0,
	\end{equation}
	then there is a real number \(\ol{\veps} = \ol{\veps}(d,\Om,Q_0)>0\)
    with the following property:
    for every \(0<\veps\leq\ol{\veps}\) and every \(f\in W^{1,p}(\Om_\veps;N)\), there exists \(\tilde{f}_\veps\in W^{1,p}(\Om;N)\) and a constant $C = C(d,p,Q_0,N)>0$ independent of $\veps$ such that
 	\begin{gather}
	\tilde{f}_\veps = f\quad \text{\ae in } \Om_\veps, \label{eq:id_ext}\\
	\norm{\tilde{f}_\veps}{L^p(\Om;N)} \le C\norm{f}{L^p(\Om_\veps;N)}, \label{eq:funct_ext}\\
	\norm{D\tilde{f}_\veps}{L^p(\Om;\R^{l\times d})} \le C\norm{Df}{L^p(\Om_\veps;\R^{l\times d})} \label{eq:grad_ext}.
	\end{gather}
\end{Theorem}

As we anticipated in the introduction,
the proof rests essentially on two results:
the extension theorem for the unconstrained case (Corollary~\ref{cor:reflect}), and
the existence of a $W^{1,p}$-retraction on the manifold $N$ (Lemma~\ref{lem:lip_retract}).
Note also that the statement of Theorem~\ref{thm:target_ext} is twofold:
on the one hand, the existence of the extension is concerned;
on the other, $L^p$-bounds for the extended function and its gradient,
respectively \eqref{eq:funct_ext} and \eqref{eq:grad_ext},
are provided.
While the construction of the extension -- hence its existence -- and
the estimate for the gradient in \eqref{eq:grad_ext}
are obtained by combining Corollary~\ref{cor:reflect} and Lemma~\ref{lem:lip_retract},
some additional work is needed to derive estimate \eqref{eq:funct_ext}.
We therefore split the proof of Theorem~\ref{thm:target_ext}
into two separate lemmas,
see Lemma~\ref{lemma:existence+grad} and Lemma~\ref{cor:Lp} below,
one dealing with the existence and the estimate on the gradient,
the other with the estimate on the function.

Before finally diving into the proof of Theorem~\ref{thm:target_ext},
to the purpose of further clarifying its content, we collect here some comments.

\begin{remark}\label{rmk:no_oper}
Theorem~\ref{thm:target_ext} gives the existence of suitable extensions. However, it does not follow from the proof (which uses a compactness argument) that an \textit{extension operator} in the spirit of \cite{AcPiMaPe92} can be defined. Since $\tilde{f}_\veps$ is obtained as the weak $W^{1,p}$-limit of a subsequence of suitable retractions, different subsequences might give different weak limits. Our construction method produces one possible extension, and the specific one we get depends on the choice made during the construction (the subsequence).
We provide a small example that roughly illustrates how the extension in Theorem~\ref{thm:target_ext} is constructed, as well as why different extensions might exist. For \(1<p<d\), let \(N\coloneqq\sph^{d-1}\). In this case, the retraction onto \(N\) in Lemma~\ref{lem:lip_retract} is the map \(P(y)\coloneqq y/|y|\) and the singular set \(X\) is a singleton given by \(X=\{0\}\). We assume that \(\Om_\veps\) is such that the origin lies in one of its open pores, i.e., there exists a radius \(\delta>0\) such that the \(d\)-dimensional ball \(B_\delta^d\) is contained in this pore. Now, let \(f:\Om_\veps\to N\) be smooth and suppose that its extension \(\sfT_\veps f\in W^{1,p}(\Om;\R^l)\) provided by Corollary~\ref{cor:reflect} takes the value zero on \(B_\delta^d\) (compare Remark~\ref{rmk:small_diff}). In this case, we cannot simply concatenate \(\sfT_\veps f\) with \(P\). However, by the arguments presented in the proof of Theorem~\ref{thm:target_ext}, for almost every sufficiently small \(h\in \R^l\) the composition \(P_h\circ \sfT_\veps f\) with the translated retraction \(P_h(y)\coloneqq P(y-h)\) is actually well-defined. For any such \(h\), the function \(\left( P_{h}|_N\right)^{-1}\circ P_{h}\circ \sfT_\veps f\) provides an admissible extension fulfilling \eqref{eq:id_ext}--\eqref{eq:grad_ext}, and Theorem~\ref{thm:target_ext} will just pick one of them. 
\end{remark}
 
\begin{remark}\label{remark:connect_assump}
Assumptions such as \eqref{eq:pi_target} have already appeared in the literature on Sobolev maps between compact Riemannian manifolds $M$ and $N$.
If $M$ has boundary $\p M$ and \(p<d\),
the $\lfloor p-1 \rfloor$-connectedness is sufficient
for $N$-valued Sobolev functions defined on $\partial M$
to admit an extension to the inside of $M$
(see \cite{HaLi87,Be14}).
A similar condition is also employed
when studying the density of smooth mappings
(cf.\ \cite{Ha09} and the references therein)
and in the statement of the Hurewicz Theorem
(cf.\ \cite[Theorem~4.32]{Ha02} or \cite[p.~57, Theorem~9.1]{Hu59}).
\end{remark}

\begin{remark}\label{rmk:lpd}
    Note that, comparing assumptions \eqref{eq:pi_target} and \eqref{eq:pi}, we see that
    Lemma~\ref{lem:lip_retract} only covers the case \(\lfloor p-1\rfloor\le l-2\). In contrast, the statement of Theorem~\ref{thm:target_ext} is formulated for \(p<d\). Thus, when \(l\le p<d\) the construction in Lemma~\ref{lem:lip_retract} does not immediately provide a retraction onto \(N\) (as well as the bound in Lemma~\ref{lem:lip_retract}\,\ref{item:projection_grad_lp}). 
    However, if for \(p\ge l\) the first \(\lfloor p-1\rfloor\) homotopy groups are trivial, this holds true in particular for the first \(l-1\) of them, and the proof of Lemma~\ref{lem:lip_retract} (see Appendix~\ref{sec:retract_lem}) shows that there exists a Lipschitz retraction \(r\in C^{0,1}(Q_R\cup \ol{N_\delta};\ol{N_\delta})\) onto the tubular neighborhood of \(N\). Setting \(P\coloneqq \xi \circ r\) for the nearest point projection \(\xi\in C^1(\ol{N_\delta};N)\), we readily obtain a Lipschitz retraction onto \(N\) that does not even possess singularities.
\end{remark}

\begin{remark}\label{remark:const}
    The construction of the extension is independent of the size of $\veps>0$. The smallness condition $\veps \leq \ol\veps$ is a consequence of the technique we use to obtain the $L^p$-bound \eqref{eq:funct_ext}, cf.\ Lemma~\ref{lemma:existence+grad} (which holds for every $\veps>0$) and Lemma~\ref{cor:Lp}.
    Actually, when \(0\in N\), estimate \eqref{eq:funct_ext} (and thus Theorem~\ref{thm:target_ext}) holds without enforcing any smallness assumption on $\veps$.
    If \(0\notin N\), the derivation of \eqref{eq:funct_ext} relies instead on an estimate of the ratio between $\mathcal{L}^d(\Omega_{0,\veps})$ and $\mathcal{L}^d(\Omega_\veps)$, hence on Lemma~\ref{lemma:estimate_measure}.
\end{remark}

\begin{remark}\label{rmk:cont}
    According to Remark~\ref{rmk:no_oper}, instead of an operator, we are describing a process or a correspondence that, for each input function $f \in W^{1,p}(\Om_\veps;N)$, produces $\tilde{f}_{\veps} \in W^{1,p}(\Om;N)$ that belongs to the set of functions that are extensions of $f$. We do not claim that this correspondence is linear. In fact, the target manifold condition itself is already non-linear in nature.
    Nevertheless, a result on the continuity of this correspondence is contained in Corollary~\ref{cor:continuity}. Note that even if the construction of an extension for a single function is not unique, the sequence of extensions behaves well in the weak topology.
\end{remark}

We now come to the case of manifold-valued Sobolev maps defined
on the perforated domain defined as in \eqref{eq:stiff_gen}.
As hinted at in Section~\ref{sec:prelim}, in comparison to Theorem~\ref{thm:target_ext},
the main difficulty in this scenario comes from the fact that intersections between the holes and the boundary \(\p \Om\) might create irregularities.
The approach employed in \cite{AcPiMaPe92} -- 
the use of partitions of unity --
is not viable in our setting, since it would result in the violation of the manifold constraint.
Therefore, in contrast to \cite{AcPiMaPe92},
we will construct extensions that are defined only on relatively compact subsets
\(\Om(\lambda)=\{x\in \Om:\ \dist(x,\p \Om)>\lambda\}\) of \(\Om\), for \(\lambda>0\).
Then, building on Theorem~\ref{thm:reflect_ext}, we are able to derive estimates on the retracted set \(\Om(\lambda)\)
as soon as \(\Om(\lambda) \subset \Om(\veps\mu)\).
Note that this approach will also enforce a smallness condition on \(\veps\),
which is however different in spirit from that of Theorem~\ref{thm:target_ext}, and it is needed also in the case $0 \in N$, cf.\ Remark~\ref{remark:const}.

\begin{corollary}\label{cor:target_ext}
    Let \(\Om\subset \R^d\), $d\geq 2$, be a bounded open set, let \(\Om_\veps^\star\) be as in \eqref{eq:stiff_gen}, and let be $\mu>0$ as in Theorem~\ref{thm:reflect_ext}.
    Under the assumptions on $p$ and \(N\) of Theorem~\ref{thm:target_ext},
    for every \(\lambda>0\), \(0<\veps\leq\lambda/\mu\), and $f \in W^{1,p}(\Om_\veps^\star;N)$, there exist $\tilde{f}_\veps^\star \in W^{1,p}(\Om(\lambda);N)$ and a constant \(C=C(d,p,E,N)>0\) independent of \(\veps\) and \(\Om\) such that
	\begin{gather}
	\tilde{f}_\veps^\star = f\quad \text{\ae in } \Om(\lambda)\cap E, \label{eq:id_ext_cor}\\
	\norm{\tilde{f}_\veps^\star}{L^p(\Om(\lambda);N)} \le C\norm{f}{L^p( \Om_\veps^\star;N)}, \label{eq:funct_ext_cor}\\
	\norm{D\tilde{f}_\veps^\star}{L^p(\Om(\lambda);\R^{l\times d})} \le C\norm{Df}{L^p(\Om_\veps^\star;\R^{l\times d})}. \label{eq:grad_ext_cor}
	\end{gather}
\end{corollary}

The lack of linearity prevents us from inferring continuity from boundedness, see Remark~\ref{rmk:cont}.
Moreover, the continuity properties do not emerge patently from the construction.
Indeed, a crucial step in the proof of Theorem~\ref{thm:target_ext} consists in choosing suitable translations of the retraction from Lemma~\ref{lem:lip_retract} and
applying them to special smooth approximations
of the extensions obtained through the operator $\sfT_\veps$ in Corollary~\ref{cor:reflect}
(see \eqref{eq:smooth_approx} and \eqref{eq:bounds_der}).
Therefore, even in the case we focus on a constant sequence in \(W^{1,p}(\Om_\veps;N)\),
we cannot be sure that we always pick the same extension along the sequence.
However, by choosing the approximations of \(\sfT_\veps f\) suitably, we are able to establish the following strong-weak continuity result (see also Remark~\ref{rmk:cont}).

\begin{corollary}\label{cor:continuity}
    Let \(\Om\), \(\Om_\veps\), $p$, and $N$ be as in Theorem~\ref{thm:target_ext}. Let $f\in W^{1,p}(\Om_\veps;N)$ and \(\sequence{f}{k}\subset W^{1,p}(\Om_\veps;N)\) be such that $f_k \to f$ in $W^{1,p}(\Om_\veps;N)$.
    Then, for every \(\veps>0\), there exists a (non-relabeled) subsequence of \(\sequence{f}{k}\) and extensions $\tilde{f}_\veps$, $(\tilde{f}^k_\veps)_k$ from Theorem~\ref{thm:target_ext} that satisfy
     \begin{equation*}
        \tilde{f}^k_\veps \rightharpoonup \tilde{f}_\veps
        \quad \text{ in } W^{1,p}(\Om;N).
    \end{equation*}
\end{corollary}

A similar result holds for the setup in Corollary~\ref{cor:target_ext} if the parameter \(\lambda\) is kept fixed.


\subsection{Proof of the existence results}\label{sec:proof_existence}

By combining Corollary~\ref{cor:reflect} and Lemma~\ref{lem:lip_retract},
we first prove the existence of extensions that satisfy \eqref{eq:id_ext} and \eqref{eq:grad_ext}. 
The $L^p$-bound in \eqref{eq:funct_ext}, which is the subject of Lemma~\ref{cor:Lp}, is then derived from geometric arguments and Lemma~\ref{lemma:estimate_measure}.

\begin{lemma}\label{lemma:existence+grad}
    Under the assumptions of Theorem~\ref{thm:target_ext}, for every $\veps >0$ and every \(f\in W^{1,p}(\Om_\veps;N)\), there exists \(\tilde{f}_\veps\in W^{1,p}(\Om;N)\) and a constant $C = C(d,p,Q_0,N)>0$ independent of $\veps$ such that \eqref{eq:id_ext} and \eqref{eq:grad_ext} hold.
\end{lemma}

\begin{proof}
    We first prove the existence of an extension, and then the $L^p$-bound \eqref{eq:grad_ext} for the gradient.

Fix $\veps>0$ and let \(f\in W^{1,p}(\Om_\veps;N)\). Applying Corollary~\ref{cor:reflect}, we can extend \(f\) to a function \(f_\veps \coloneqq\sfT_\veps f\in W^{1,p}(\Om;\R^l)\). 
Note that, at this point, the values of \(f_\veps\) are not necessarily contained in \(N\). Therefore, one would like to compose $f_\veps$ with a suitable retraction onto $N$. When applying Lemma~\ref{lem:lip_retract} for this purpose, we have to make sure that $f_\veps$ does not take values in the singularity set \(X\) on a subset of positive Lebesgue measure. Unfortunately, this possibility cannot be excluded a priori. However, inspired by \cite[Theorem~6.2]{HaLi87}, we avoid this problem by working with smooth approximations of $f_\veps$ and possibly perturbing the retraction from Lemma~\ref{lem:lip_retract} by a small translation.

The Meyers-Serrin Theorem (see, e.g., \cite[Theorem~11.24]{Leoni17}) guarantees the existence of a sequence \((f_\veps^k)_k\subset C^\infty(\Om;\R^l) \cap W^{1,p}(\Om;\R^l)\) such that
\begin{equation}\label{eq:smooth_approx}
	f^k_\veps \to f_\veps \ \mbox{ in } W^{1,p}(\Om;\R^l).
\end{equation}
	
Recall that \(f(x)\in N\) for \ae \(x\in \Om_\veps\) and $N$ is compact by assumption. This, together with Remark~\ref{rmk:small_diff}, implies that the function \(f_\veps\) only takes values in some compact set, not depending on the special choice of \(f\). Thus, we can pick the approximating sequence \((f^k_\veps)_k\) in such a way that the image of each of the functions is contained in some \(l\)-dimensional cube \(Q_{\hat{R}} = (-\hat{R},\hat{R})^l\), for some \(\hat{R}>0\) independent of \(k\in\N\). In fact, the functions constructed in the Meyers-Serrin approximation are exactly of this type since they are obtained by a convolution and partition procedure.
	
Let now $N_\delta$ be the tubular neighborhood of $N$ defined in Section~\ref{sec:prelim}, and let $P : Q_R \setminus X \to N$ be given by Lemma~\ref{lem:lip_retract} with $j = \lfloor p-1\rfloor$ (see also Remark~\ref{rmk:lpd}). Note that we can always assume $R > \hat{R}$. Indeed, by enlarging the $l$-dimensional cube $Q_R$ containing $N_\delta$ we will possibly obtain a different singular set $X$ and a different retraction $P$, but still such that assertions \ref{item:retraction}--\ref{item:projection_grad_lp} hold.
	
We set $\sigma \coloneqq \min\{R-\hat{R},\delta\}$ and, for $h \in B_\sigma \coloneqq \{y \in \R^l : |y| < \sigma\}$, we define the translated retractions 
\begin{align}
\begin{split}\label{eq:translated_proj}
		P_h :\, &(Q_R+h) \setminus (X+h) \to N, \\
		&P_h(y) \coloneqq P(y-h). 
\end{split}
\end{align}
These are well-defined, since $y-h \in Q_R\setminus X$ whenever $y\in (Q_R+h)\setminus (X+h)$. Moreover, for all $h \in B_\sigma$ we have
\begin{equation}\label{eq:inclusions}
		N-h\subset N_\delta,
\end{equation}
so that $N \subset (Q_R+h) \setminus (X+h)$. We would now like to compose $f^k_\veps$ with $P_h$. Note that for all $h \in B_\sigma$ we have $Q_{\hat{R}} \subset Q_R+h$, so that the image of $f^k_\veps$ actually belongs to $Q_R+h$. We need to show that it falls outside the perturbed singular sets $X+h$ for \ae $h \in B_\sigma$, $\mathcal{L}^d$-almost everywhere.
By construction, the set $X$ is contained in a finite union of $l-\lfloor p \rfloor +1$-dimensional planes. Then, as a consequence of Sard's theorem (see \cite[Lemma~2.3]{BPV14}), for any \(k\in\N\) and \ae $h \in \R^l$ the preimage of $X+h$ under $f^k_\veps$ is a finite union of smooth submanifolds of dimension $d-\lfloor p+1 \rfloor$, hence $\mathcal{L}^d$-negligible.
Consequently, we obtain that for \ae $h \in B_\sigma$ the compositions $P_h \circ f^k_\veps$ are well-defined almost everywhere and weakly differentiable.
	
We now prove that there exists some \(h_k \in B_\sigma\) such that \(P_{h_k} \circ f^k_\veps \in W^{1,p}(\Omega;N)\). To that end, first observe that
\begin{equation}\label{eq:comp_lp}
	\int_{\Om} |P_h \circ f^k_\veps(x)|^p {\dl x} 
	\le C\mathcal{L}^d(\Om)
	<\infty
	\quad \mbox{for all } h \in B_\sigma,
\end{equation}
with some constant $C>0$ independent of $\veps$ and $l$. This comes from the fact that \(P_h \circ f^k_\veps\in N\) \ae by construction, and \(N\) as well as \(\Om\) are bounded by assumption. Second, we compute with Fubini's theorem, the chain rule, and \ref{item:projection_grad_lp} of Lemma~\ref{lem:lip_retract} that 
\begin{align}
\begin{split}\label{eq:comp_der_lp}
	\int_{B_\sigma}\int_{\Om} |D(P_h \circ f^k_\veps)(x)|^p {\dl x} {\dl h}
	&\le \int_{\Om} |Df^k_\veps(x)|^p \int_{B_\sigma} |DP_h(f^k_\veps(x))|^p {\dl h} {\dl x}\\
	&= \int_{\Om} |Df^k_\veps(x)|^p \int_{B_\sigma} |DP(f^k_\veps(x)-h)|^p {\dl h} {\dl x}\\
	&\le \int_{\Om} |Df^k_\veps(x)|^p \int_{Q_R} |DP(y)|^p {\dl y} {\dl x}\\
	&\le
	C\norm{Df^k_\veps}{L^p(\Om;\R^{l\times d})}^p.
\end{split}
\end{align}
In particular, from the mean value theorem, there has to exist some \(h_k \in B_\sigma\) such that
\begin{equation}\label{eq:bounds_der}
	\norm{D(P_{h_k} \circ f^k_\veps)}{L^p(\Om;\R^{l\times d})}
	\le \frac{C}{\big(\mathcal{L}^d(B_\sigma)\big)^{1/p}} \norm{Df^k_\veps}{L^p(\Om;\R^{l\times d})}.
\end{equation}
Note that, together with \eqref{eq:comp_lp}, we infer that indeed \(P_{h_k} \circ f^k_\veps\in W^{1,p}(\Om;N)\) for every $k\in\N$.
	
The translated retractions have solved the problem that the functions \(P\circ {f^k_\veps}\) might be ill-defined. At the same time, however, they may contradict the requirement that the extension is equal to the original function on \(\Om_\veps\). In fact, there could exists \(x\in\Om\) such that \(P_{h_k}\circ f^k_\veps(x)\ne f^k_\veps(x)\), even though \(f^k_\veps(x)\in N\). This is because \(P_{h_k}|_N\) is not equal to the identity. It turns out that with the additional composition
\begin{equation}\label{eq:final_pro}
	\left( P_{h_k}|_N\right)^{-1}\circ P_{h_k}\circ f^k_\veps \ : \ \Om\to N
\end{equation}
the issue is solved, as indeed
\begin{equation}
	\left( \left( P_{h_k}|_N\right)^{-1}\circ P_{h_k}\circ f^k_\veps\right)(x)
	=f^k_\veps(x)
	\quad \text{for all } x\in\Om \text{ such that } f^k_\veps(x)\in N.
\end{equation}
Note that \(P_h|_N(\,\cdot\,) = P|_{N-h}(\,\cdot\, - h)\), where, thanks to \eqref{eq:inclusions}, the function \(P|_{N-h} : N-h\to N\) is just the nearest point projection for any \(h \in B_\sigma\). Thus, from the assumption that \(N\) is \(C^2\) and the regularity of the nearest point projection (compare Lemma~\ref{lem:nearest_pp}), we can infer that \(P_h|_N \in C^1(N;N)\). Furthermore, switching to local coordinates, it is readily checked that \(P_h|_N\) is invertible with inverse function \(\left( P_h|_N\right)^{-1}=P_{-h}|_N \in C^1(N;N)\). In particular,
\(\left( P_h|_N\right)^{-1}\) is Lipschitz on \(N\), and
\begin{equation}\label{eq:inverse_der_lp}
	\Lambda \coloneqq \sup_{h\in B_\sigma} \textup{Lip}\Big(\left( P_h|_N\right) ^{-1}\Big) < \infty
\end{equation} 
is a finite number depending only on $N$. We conclude, together with an argument similar to \eqref{eq:comp_lp}, 
that \(\left( P_{h_k}|_N\right)^{-1}\circ P_{h_k}\circ f^k_\veps\in W^{1,p}(\Om;N)\) for all \(k\in \N\). We are now left to prove that the resulting sequence converges in an appropriate way.
		
It follows immediately from our calculations that \(((P_{h_k}|_N)^{-1}\circ P_{h_k} \circ f^k_\veps)_k\) is uniformly bounded in \(W^{1,p}(\Om;N)\). Hence, by compactness, there exists some \(\tilde{f}_\veps\in W^{1,p}(\Om;\R^l)\) such that
\begin{equation}\label{eq:weak_conv}
		(P_{h_k}|_N)^{-1}\circ P_{h_k} \circ f^k_\veps
		\rightharpoonup \tilde{f}_\veps \ \mbox{ in } W^{1,p}(\Om;\R^l).
\end{equation}
By choosing a pointwise \ae converging subsequence, we can infer that \(\tilde{f}_\veps \in W^{1,p}(\Om;N)\). Additionally, since \((h_k)_k\subset B_\sigma\), by the Bolzano-Weierstrass theorem there exists \(\tilde{h} \in B_\sigma\) such that \(h_k\to \tilde{h}\) (up to a subsequence). Now, note that \(\cup_{k} (f^k_\veps)^{-1}(X+h_k)\) is a \(\mathcal{L}^d\)-null-set. Thus, we may assume that \ae \(x\in \Om \setminus (f_\veps)^{-1}(X+\tilde{h})\) is not contained in any of the sets \((f^k_\veps)^{-1}(X+h_k)\). This implies that, for such an $x$, the expressions \((P_{{h_k}}|_N)^{-1}\circ P_{{h_k}} \circ f^k_\veps(x)\) and \((P_{\tilde{h}}|_N)^{-1}\circ P_{\tilde{h}} \circ f_\veps(x)\) are well-defined.
Consequently,
since \(f^k_\veps \to f_\veps\) \ae as \(k\to \infty\), the continuity of $P$ yields
\begin{equation*}
\begin{aligned}
		&\lim_{k \to \infty} \Big((P_{h_k}|_N)^{-1} \circ P_{h_k} \circ f^k_\veps(x)\Big) = \lim_{k \to \infty} P\big(P(f^k_\veps(x) - h_k) + h_k\big) \\
		&= P\big(P(f_\veps(x) - \tilde{h}) + \tilde{h}\big) = (P_{\tilde{h}}|_N)^{-1} \circ P_{\tilde{h}} \circ f_\veps(x)
		\quad \text{ for \ae } x \in \Omega \setminus (f_\veps)^{-1}(X+\tilde{h}).
\end{aligned}
\end{equation*}
It follows now by uniqueness of the weak \(W^{1,p}\)-limit that for \ae \(x\in \Om\setminus (f_\veps)^{-1}(X+\tilde{h})\) the limit function \(\tilde{f}_\veps\) is given by
\begin{equation}\label{eq:def_ext}
		\tilde{f}_\veps(x)
		= (P_{\tilde{h}}|_N)^{-1}\circ P_{\tilde{h}} \circ f_\veps(x).
\end{equation}
Since \(f_\veps=f\) \ae on \(\Om_\veps\) by \eqref{eq:ex_op_simplified_prop1}, \(f(x)\in N\) for \ae \(x\in \Om_\veps\), and \(N\cap (X+\tilde{h}) = \emptyset\) from \ref{item:inclusion} of Lemma~\ref{lem:lip_retract} and \eqref{eq:inclusions}, for \ae $x\in \Om_\veps$ we get 
\begin{align*}
	\tilde{f}_\veps(x)
	= (P_{\tilde{h}}|_N)^{-1}\circ P_{\tilde{h}} \circ f_\veps(x)
	= (P_{\tilde{h}}|_N)^{-1}\circ P_{\tilde{h}} \circ f(x)
	= f(x),
\end{align*}
so that \eqref{eq:id_ext} holds.

Estimate \eqref{eq:grad_ext} for the gradient follows directly from \eqref{eq:bounds_der} and \eqref{eq:inverse_der_lp} together with the convergences \eqref{eq:smooth_approx} and \eqref{eq:weak_conv}, concatenated with the corresponding $L^p$-estimate for $D f_\veps=D(\sfT_\veps f)$ in \eqref{eq:ex_op_simplified_prop3}.
\end{proof}

As a byproduct of the preceding proof,
we obtain the following pointwise characterization of the extension
(see \eqref{eq:def_ext}).

\begin{corollary}\label{cor:aux_comp}
Under the assumptions of Theorem~\ref{thm:reflect_ext}, there exists \(0 < \tilde{h} < \delta\) such that for \ae \(x\in \Omega\setminus (\sfT_\veps f)^{-1}(X+\tilde{h})\) the extension \(\tilde{f}_\veps\) is given by
\begin{equation}\label{eq:def_ext1}
    \tilde{f}_\veps(x)
    = (P_{\tilde{h}}|_N)^{-1}\circ P_{\tilde{h}} \circ \sfT_\veps f(x),
\end{equation}
where \(\sfT_\veps:W^{1,p}(\Om_\veps;\R^l) \to W^{1,p}(\Om;\R^l)\) is the extension operator from Corollary~\ref{cor:reflect}, and $P_{\tilde{h}}$ is defined as \(P_{\tilde{h}}(y) = P(y-\tilde{h})\) for \(P\) the retraction from Lemma~\ref{lem:lip_retract}.
\end{corollary}

We now turn to the $L^p$-estimate in \eqref{eq:funct_ext}.
\begin{lemma}\label{cor:Lp}
    Under the assumptions of Theorem~\ref{thm:target_ext},
    there is a real number \(\ol{\veps} = \ol{\veps}(d,\Om,Q_0)>0\)
    with the following property: for every \(f\in W^{1,p}(\Om_\veps;N)\), the extension $\tilde{f}_\veps \in W^{1,p}(\Om;N)$ from Lemma~\ref{lemma:existence+grad}
    satisfies \eqref{eq:funct_ext} 
    for every \(\veps\leq\ol{\veps}\) and some constant $C = C(d,p,Q_0,N)>0$ independent of $\veps$.
\end{lemma}

\begin{proof}
For $\veps >0$ and $f \in W^{1,p}(\Om_\veps;N)$, let $\tilde{f}_\veps \in W^{1,p}(\Om;N)$ be
the extension from Lemma~\ref{lemma:existence+grad}.
The derivation
of the \(L^p\)-estimate in \eqref{eq:funct_ext}
relies on the properties of the retraction on $N$ from Lemma~\ref{lem:lip_retract}.
Precisely, let $P_{\tilde{h}}$ be as in~\eqref{eq:def_ext1}.
Then, the following hold:
\begin{enumerate}[nolistsep]
    \item\label{item:lipschizt_prop}
    the function 
    \begin{equation}\label{eq:comp_proj}
         \hat{P}	:\, N_\delta \to N,
        \quad 
         \hat{P}\coloneqq\left( P_{\tilde{h}}|_N\right)^{-1} \circ P_{\tilde{h}}|_{N_\delta}
    \end{equation}
    is globally Lipschitz continuous on $\overline{N_\delta}$ by Lemma~\ref{lem:nearest_pp},
    \item\label{item:identity_prop}
    and $\hat{P}$ is the identity on $N$.
\end{enumerate}

We define
\begin{equation}
   K_{\veps,\delta}
   \coloneqq \{x\in \Om_{0,\veps} : \, \sfT_\veps f(x)\in N_\delta\},
\end{equation}
where, as before, $\Om_{0,\veps}$ is as in \eqref{eq:second_setup2} and \(\sfT_\veps\) is the extension operator from Corollary~\ref{cor:reflect}. The set \(K_{\veps,\delta}\) is measurable by construction.

Let us fix $y\in N$.
On the one hand, by \eqref{eq:def_ext} and by properties \eqref{item:lipschizt_prop} and \eqref{item:identity_prop} above, we have
\begin{equation}\label{eq:LipP}
    \int_{K_{\veps,\delta}}| \tilde{f}_\veps - y |^p {\dl x}
    = \int_{K_{\veps,\delta}}| \hat{P}\circ \sfT_\veps f - \hat{P}(y) |^p {\dl x}
    \le \left(\textup{Lip}\hat{P}\right)^p \int_{K_{\veps,\delta}}|\sfT_\veps f - y |^p {\dl x}.
\end{equation}
On the other hand, letting \(K_{\veps,\delta}^c\) be the complement of \(K_{\veps,\delta}\) in \(\Om_{0,\veps}\),
that is, \(K_{\veps,\delta}^c \coloneqq \{x\in \Om_{0,\veps}:\, \sfT_\veps f(x)\notin N_\delta\}\), and
setting \(\Gamma \coloneqq \max_{y\in N} | y |\),
we also find
\begin{equation}\label{eq:Ncompact}
    \int_{K^c_{\veps,\delta}}| \tilde{f}_\veps - y |^p {\dl x}
    \le \left(\frac{2\Gamma}{\delta}\right)^p \int_{K^c_{\veps,\delta}}\delta^p {\dl x}
    \le \left(\frac{2\Gamma}{\delta}\right)^p \int_{K^c_{\veps,\delta}} | \sfT_\veps f - y |^p {\dl x},
\end{equation}
where the second inequality follows from the definition of $K^c_{\veps,\delta}$.
By summing up \eqref{eq:LipP} and \eqref{eq:Ncompact},
for a suitable constant \(C=C(N)\ge0\) we obtain
\begin{equation*}
        \int_{\Om_{0,\veps}} | \tilde{f}_\veps - y |^p {\dl x}
        \le C \int_{\Om_{0,\veps}}| \sfT_\veps f - y |^p {\dl x},
\end{equation*}
and we thereby derive the following chain of estimates: 
\begin{align*}
    \norm{\tilde{f}_\veps}{L^p(\Om_{0,\veps};N)} & \le
    \norm{\tilde{f}_\veps - y }{L^p(\Om_{0,\veps};N)} + \Big(\mathcal{L}^d(\Om_{0,\veps})\Big)^{1/p}|y| \\
    & \le C \norm{\sfT_\veps f - y }{L^p(\Om_{0,\veps};N)} + \Big(\mathcal{L}^d(\Om_{0,\veps})\Big)^{1/p}|y| \\
    & \le C \norm{\sfT_\veps f }{L^p(\Om_{0,\veps};N)} + 2 \Big(\mathcal{L}^d(\Om_{0,\veps})\Big)^{1/p}|y|.
\end{align*}
Since $y$ is an arbitrary point on the compact manifold $N$,
we can optimize the last quantity so as to get
\begin{align*}
    \norm{\tilde{f}_\veps}{L^p(\Om_{0,\veps};N)} \le C \norm{\sfT_\veps f }{L^p(\Om_{0,\veps};N)} + 2 \gamma \Big(\mathcal{L}^d(\Om_{0,\veps})\Big)^{1/p},
\end{align*}
where \(\gamma \coloneqq \min_{y\in N} | y |\).
Recalling that \(\tilde{f}_\veps = \sfT_\veps f\) on \(\Om_\veps\),
from the previous estimate we obtain
\begin{equation}\label{eq:gamma0}
    \norm{\tilde{f}_\veps}{L^p(\Om;N)} \le C \norm{\sfT_\veps f }{L^p(\Om;N)} + 2 \gamma \Big(\mathcal{L}^d(\Om_{0,\veps})\Big)^{1/p},
\end{equation}
whence, by the property of the operator $\sfT_\veps$ in \eqref{eq:ex_op_simplified_prop2},
we infer
\begin{equation}\label{eq:gamma1}
    \norm{\tilde{f}_\veps}{L^p(\Om;N)} \le C \norm{ f }{L^p(\Om_\veps;N)} + 2\gamma \Big(\mathcal{L}^d(\Om_{0,\veps})\Big)^{1/p}
\end{equation}
for a constant $C>0$ independent of $\veps$ and $\Omega$.
In particular, in the case $0\in N$, $\gamma = 0$, and we immediately obtain \eqref{eq:funct_ext}.

If \(0\notin N\), instead, a further step is needed. We can write
\begin{equation}\label{eq:gamma2}
    \gamma\Big(\mathcal{L}^d(\Om_{0,\veps})\Big)^{1/p}
    = \left(\frac{\mathcal{L}^d(\Om_{0,\veps})}{\mathcal{L}^d(\Om_\veps)}\right)^{1/p}\left(\int_{\Om_\veps}\gamma \dl x\right)^{1/p}
    \le \left(\frac{\mathcal{L}^d(\Om_{0,\veps})}{\mathcal{L}^d(\Om_\veps)}\right)^{1/p}\norm{ f }{L^p(\Om_\veps;N)}.
\end{equation}
Then, thanks to Lemma~\ref{lemma:estimate_measure} we deduce that
there exists $\ol{\veps} = \ol{\veps}(d,\Om,Q_0)>0$ such that
    \[
    \gamma\Big(\mathcal{L}^d(\Om_{0,\veps})\Big)^{1/p} \leq C \norm{ f }{L^p(\Om_\veps;N)}
    \qquad\text{for all } \veps \leq \ol{\veps}, 
    \]
where $C = C(p,Q_0)$ is a suitable constant.
By plugging the last inequality into \eqref{eq:gamma1},
we achieve the conclusion.
\end{proof}

At this stage, the proof of Theorem~\ref{thm:target_ext} readily follows.

\begin{proof}[Proof of Theorem~\ref{thm:target_ext}]
    The existence of an extension was settled in Lemma~\ref{lemma:existence+grad}.
    The same lemma shows that estimate~\eqref{eq:grad_ext} holds.
    Thus, the desired result is achieved by applying Lemma~\ref{cor:Lp},
    since it establishes estimate~\eqref{eq:funct_ext} too.
\end{proof}

To conclude this subsection, let us now turn to the proof of Corollary~\ref{cor:target_ext}, that is, the construction of suitable extensions for manifold-valued Sobolev functions when the perforated domain is defined as in \eqref{eq:stiff_gen}.

\begin{proof}[Proof of Corollary~\ref{cor:target_ext}]
Let us fix $\lambda > 0$ and let $\sfT_\veps^\star$ be the extension operator from Theorem~\ref{thm:reflect_ext}.
By the same steps as in the proof of Theorem~\ref{thm:target_ext},
but now applied to the restriction to \(\Om(\lambda)\) of the extended maps \(\sfT_\veps^\star f\)
-- approximating \(\sfT_\veps^\star f\) by smooth functions, extending these, composing the extension with suitable retractions, and then taking the limit --, we can construct a function \(\tilde{f}_\veps^\star \in W^{1,p}(\Om(\lambda);N)\) satisfying \eqref{eq:id_ext_cor}.

It remains to verify the claimed \(L^p\)-estimates.
As for the gradient, we note that by Theorem~\ref{thm:reflect_ext}
we have an estimate for \(D(\sfT_\veps^\star f)\) on the retracted set \(\Om(\veps \mu)\).
However, under the assumption \(\veps\le \lambda/\mu\),
it holds \(\Om(\lambda)\subset \Om(\veps \mu)\) and therefore,
in view of~\eqref{eq:ex_op_prop3} and similarly to the proof of Theorem~\ref{thm:target_ext},
we infer \eqref{eq:grad_ext_cor}.

Finally, let us turn to the proof of \eqref{eq:funct_ext_cor}.
An analogue to \eqref{eq:gamma0} on \(\Om(\lambda)\) can be derived, namely, 
\begin{equation*}
        \norm{\tilde{f}_\veps^\star}{L^p(\Om(\lambda);N)} \le C \norm{\sfT_\veps^\star f }{L^p(\Om(\lambda);N)} + 2 \gamma \Big(\mathcal{L}^d\big(\Om_{0,\veps}^\star\cap\Om(\lambda)\big)\Big)^{1/p}
\end{equation*}
with $\gamma = \min_{y \in N}|y|$.
Then, again if \(\veps\le \lambda/\mu\),
from the previous inequality and the \(L^p\)-estimate \eqref{eq:ex_op_prop2} for \(\sfT_\veps^\star f\),
we obtain that 
\begin{equation}\label{eq:gamma4}
    \norm{\tilde{f}_\veps^\star}{L^p(\Om(\lambda);N)}
        \le C \norm{ f }{L^p(\Om_\veps^\star;N)} + 2\gamma\Big(\mathcal{L}^d\big(\Om_{0,\veps}^\star\cap\Om(\veps\mu)\big)\Big)^{1/p},
\end{equation}
for some constant \(C>0\) not depending on $\veps$ and $\Omega$.

Now, suppose that $\mathcal{L}^d(\Om_\veps^\star)>0$
(the other case leading to a trivial statement).
Then, similarly to \eqref{eq:gamma2}, we observe that
\begin{equation}\label{eq:gamma5}
    \gamma\Big(\mathcal{L}^d\big(\Om_{0,\veps}^\star\cap\Om(\veps\mu)\big)\Big)^{1/p}
    \le \left(\frac{\mathcal{L}^d\big(\Om_{0,\veps}^\star\cap\Om(\veps\mu)\big)}{\mathcal{L}^d(\Om_\veps^\star)}\right)^{1/p}\norm{ f }{L^p(\Om_\veps^\star;N)}.
\end{equation}
To bound the right-hand side from above,
we observe that from the definition of \(\Om_\veps^\star\) and the inclusion \(Z_\veps^{\star,\textup{int}}\subset Z_\veps^\star\) it follows 
\begin{equation*}
    \mathcal{L}^d(\Om_\veps^\star)
	\ge \sum_{Z_\veps^{\star,\textup{int}}} \veps^d \mathcal{L}^d(Q_1).
\end{equation*}
This, combined with inequality \eqref{eq:gen_soft_est}, yields
\begin{equation*}
    \frac{\mathcal{L}^d\big(\Om_{0,\veps}^\star\cap\Om(\veps \mu)\big)}{\mathcal{L}^d(\Om_\veps^\star)} 
    \le \frac{\sum_{Z_\veps^{\star,\textup{int}}} \veps^d \mathcal{L}^d(Q_0)}{\sum_{Z_\veps^{\star,\textup{int}}} \veps^d \mathcal{L}^d(Q_1)} 
    = \frac{\mathcal{L}^d(Q_0)}{\mathcal{L}^d(Q_1)}.
\end{equation*}
By inserting the last estimate into \eqref{eq:gamma5} and
by combining it with \eqref{eq:gamma4},
we obtain \eqref{eq:funct_ext_cor}.
\end{proof}


\subsection{Proof of the strong-weak continuity}
We devote this subsection to the proof of Corollary~\ref{cor:continuity}.
As we already pointed out, given a strongly converging sequence \(\sequence{f}{k}\) with limit \(f\), in order to achieve the weak convergence of the extensions \(\tilde{f}_\veps^k\) to an extension $\tilde{f}_\veps$ of $f$,
one has to select the approximating sequence for \(\sfT_\veps f\) in the proof of Lemma~\ref{lemma:existence+grad} carefully.
In doing so, we will also automatically solve the problem of picking suitable translations of the retraction in the construction of \(\tilde{f}_\veps\).

\begin{proof}[Proof of Corollary~\ref{cor:continuity}]
Let \(\sequence{f}{k}\subset W^{1,p}(\Om_\veps;N)\) be a sequence converging in \(W^{1,p}(\Om_\veps;N)\) to some limit function \(f\). To simplify the notation, let us first assume that \(\veps=1\). The proof for a general \(\veps >0\) can always be reduced to this case by a simple scaling argument similar to the one in \cite{AcPiMaPe92}, as we will show below. We will write \(\tilde{f}\) instead of \(\tilde{f}_\veps\) for the extension from Theorem~\ref{thm:target_ext}, and $\sfT$ instead of \(\sfT_\veps\) for the operator in Corollary~\ref{cor:reflect}.

As in the proof of Theorem~\ref{thm:target_ext}, for every fixed \(k\), let \((f_j^k)_j\subset C^\infty(\Om;\R^l) \cap W^{1,p}(\Om;\R^l)\) be a smooth sequence such that \(f_j^k\to \sfT f_k\) in \(W^{1,p}(\Om;\R^l)\) for \(j\to \infty\). By the continuity of the operator \(\sfT : W^{1,p}(\Om_1;\R^l) \to W^{1,p}(\Om;\R^l)\), we know that \(\sfT f_k\to \sfT f\) in \(W^{1,p}(\Om;\R^l)\) for \(k\to \infty\). Hence,
    \begin{equation*}
        \lim_{k\to\infty}\lim_{j\to\infty} \norm{f_j^k-\sfT f}{W^{1,p}(\Om;\R^l)}
        =0.
    \end{equation*}
    At the same time, for fixed \(k\), without loss of generality we can assume that, upon taking a subsequence, \(\left( P_{h_j}|_N\right)^{-1}\circ P_{h_j}\circ f^k_j\rightharpoonup \tilde{f}_k\) in \(W^{1,p}(\Om;N)\) for \(j\to \infty\), where \(\tilde{f}_k\) are the extensions of $f_k$ provided by Theorem~\ref{thm:target_ext}.\\
    Therefore, by the metrizability of weak convergence on bounded sets of \(W^{1,p}(\Om;\R^l)\), there exists a subsequence \(j(k)\to \infty\) such that \(f^k_{j(k)}\to \sfT f\) in \(W^{1,p}(\Om;\R^l)\) and 
    \begin{equation*}
        \left( P_{h_{j(k)}}|_N\right)^{-1}\circ P_{h_{j(k)}}\circ f^k_{j(k)}-\tilde{f}_k \rightharpoonup 0 \ \mbox{ in } W^{1,p}(\Om;N).
    \end{equation*}
    Since 
    \begin{equation*}
        \left( P_{h_{j(k)}}|_N\right)^{-1}\circ P_{h_{j(k)}}\circ f^k_{j(k)}
        \in W^{1,p}(\Om;N)
    \end{equation*}
    and the sequence is bounded in that space, it converges, up to another subsequence, weakly in \(W^{1,p}\) to some function \(\tilde{f}\in W^{1,p}(\Om;N)\). Then, owing to the proof of Theorem~\ref{thm:target_ext}, the desired extension of $f$ can be chosen as the function \(\tilde{f}\). This proves the statement in the case $\veps=1$.

For an arbitrary \(\veps>0\), we combine the above steps with the scaling procedure in \cite{AcPiMaPe92}. We define the invertible affine map \(\sigma_\veps:\R^d\to \R^d\) as \(\sigma_\veps(x)\coloneqq \veps x\) for $x \in \R^d$. This allows us to find extensions \(\tilde{f}_\veps^k= (\widetilde{f_k\circ \sigma_{\veps}})\circ \sigma_{1/\veps}\) and \(\tilde{f}_\veps= (\widetilde{f\circ \sigma_\veps})\circ \sigma_{1/\veps}\) such that \(\widetilde{f_k\circ \sigma_{\veps}}\rightharpoonup \widetilde{f\circ \sigma_\veps}\) in \(W^{1,p}((1/\veps) \Om;N)\). This proves Corollary~\ref{cor:continuity}.
\end{proof}


\subsection{Application: a quick approach to derive effective models in micromagnetics}\label{sec:ex}

As we touched upon in the introduction,
our motivation for studying extensions as in Theorem~\ref{thm:target_ext} roots in homogenization problems involving domains with microperforations, such as \(\Om_\veps\) above.
The goal of homogenization consists in deriving an effective, macroscopic model in the limit \(\veps\to 0\) by starting from the descriptions at the microscopic level. If at this scale the state of the system is encoded by a function \(f_\veps\) -- the minimizer of a suitable variational problem or the solution to a system of partial differential equations --, then the homogenization involves the determination of the limiting behavior of the family \(\sequence{f}{\veps}\) with respect to a suitable notion of convergence.
As we have already mentioned, the hurdle posed by the microperforations can be dealt with by constructing suitable extensions. However, any manifold constraint that is prescribed at the \(\veps\)-level is in general lost in this process (and {\it a fortiori} in the limit). 
In some specific scenarios, such a problem is actually solvable, for example by some {\it a posteriori} argument, see Remark \ref{rmk:2scale-approach} below.
Here, we propose a straightforward approach based on the manifold-constrained extensions from Theorem~\ref{thm:target_ext}.
As a proof of concept, we discuss a variational toy problem inspired by the theory of micromagnetics.

It is experimentally validated that,
below a certain temperature threshold known as the \emph{Curie temperature},
the vector field \(M\) describing the magnetization of ferromagnetic materials is subject to a {\it saturation} (or {\it Heisenberg}) constraint that fixes its norm.
Accordingly, if the sample sits in $\Om\subset\R^3$, the magnetization is usually expressed as \(M=M_s m\),
where \(m\in H^1(\Om;\sph^2)=W^{1,2}(\Om;\sph^2)\) (cf.\ \cite{Br63, HuSch08}) is a rescaled vector field.
The well-established theory of micromagnetics, proposed by Landau and Lifshitz \cite{LaLi35}, states
that stationary magnetic configurations on a perforated material sample minimize an energy functional made up of several contributions,
each of them accounting for distinct phenomena.
As a simplified model (cf.~\cite{AlDF15}),
for \(\veps >0\) we consider here the functional
\begin{equation}\label{eq:toy_funcs}
		\mathcal{I}_\veps(m)
		\coloneqq \int_{\Om_\veps} g\left(
            m, Dm\right) {\dl x},
\end{equation}
where \(m \in H^1(\Om_\veps;\sph^2)\) and the energy density \(g:\R^3\times \R^{3\times 3}\to \R\):
\begin{enumerate}[nolistsep]
    \item is continuous;
    \item satisfies \begin{equation}\label{eq:equi_coercivity}
    c |A|^2 - C
    \le g(\nu,A)
    \quad \forall (\nu,A)\in \R^3\times \R^{3\times 3}
    \end{equation}
    for some \(c,C>0\).
\end{enumerate}

Note that on the left-hand side of \eqref{eq:equi_coercivity} no term including \(|m|\) appears as a consequence of the constraint \(|m|=1\) \ae on \(\Om\).

In the present setup, \(d=3\), \(p=2\), and \(\pi_0(\sph^2)=\pi_1(\sph^2)=0\),
and the problem falls therefore into the scope of Theorem~\ref{thm:target_ext}.
By applying it, we obtain the following result for the family of micromagnetic energies in \eqref{eq:toy_funcs}.
\begin{lemma}\label{lem:micromag}
    Let the functionals $(\mathcal{I}_\veps)_\veps$ be defined as in \eqref{eq:toy_funcs}, where \(g:\R^3\times \R^{3\times 3}\to \R\) is a continuous function satisfying \eqref{eq:equi_coercivity}. If the sequence $(m_\veps)_\veps$ is such that
    \begin{equation}\label{eq:bouIeps}
        \liminf_{\veps\to0} \mathcal{I}_\veps(m_\veps)<+\infty,
    \end{equation}
    then there exist a subsequence (that we do not relabel) and $m \in W^{1,2}(\Om;\sph^2)$ such that
    $$
    \tilde{m}_\veps \rightharpoonup m \ \text{ weakly in } W^{1,2}(\Om;\sph^2),
    $$
    where $\tilde{m}_\veps$ is the extension of $m_\veps$ constructed in Theorem~\ref{thm:target_ext}.
\end{lemma}

\begin{proof}
    Thanks to \eqref{eq:bouIeps},
    it follows from estimate \eqref{eq:grad_ext} in Theorem~\ref{thm:target_ext}, the coercivity assumption \eqref{eq:equi_coercivity}, and the Banach-Alaoglu theorem that the family \((\tilde{m}_\veps)_\veps \in W^{1,2}(\Om;\mathbb{S}^{2})\) admits a subsequence
    that converges weakly to some limit function \(m\in W^{1,2}(\Om;\R^3) \).
    By compact embedding, we might extract a further subsequence that converges strongly in \(L^2(\Om;\R^3)\) and even pointwise a.\,e.\ in \(\Om\) to $m$.
    We thereby infer that there holds indeed \(m\in W^{1,2}(\Om;\mathbb{S}^{2})\). 
\end{proof}

\begin{remark}\label{rmk:2scale-approach}
    The conclusion of Lemma~\ref{lem:micromag} can also be obtained using strategies different from ours, for example, by two-scale methods \cite{All92},
    as in \cite{ChoMouTi18}.
    Nevertheless, the use of Theorem~\ref{thm:target_ext} makes the proof shorter and more direct:
    indeed, if the extensions from \cite{AcPiMaPe92} are employed,
    one needs to recover the constraint \(m\in W^{1,p}(\Om;\mathbb{S}^{2})\) for the liming energy functional
    in an {\it ad hoc} separate step,
    cf.\ \cite[Lemma~4.1]{ChoMouTi18}.\newline
   Note also that the strategy in \cite{ChoMouTi18} may be followed for arbitrary submanifolds $N$ of \(\R^l\) that coincide with the level set of some continuous function $h : \R^l \to \R^{l-n}$.
   Our approach, instead, provides a straightforward way to preserve the manifold constraint in the limiting model for every $\lfloor p-1 \rfloor$-connected, compact target manifold without boundary.
\end{remark}

\begin{remark}
    Once Lemma \ref{lem:micromag} is on hand,
    it can be inferred that
    the domain of the energy functional that provides the homogenized, macroscopic description of a microperforated micromagnetic body is \(W^{1,2}(\Om;\mathbb{S}^{2})\).
    In the current variational setting,
    this follows by standard $\Gamma$-convergence techniques,
    see, e.g., \cite[Chapter 19]{BrDe98}.
    If the limit of \((\mathcal{I}_\veps)_\veps\) in the sense of $\Gamma$-convergence is computed with respect to the weak topology of $W^{1,2}(\Omega;\sph^2)$ (equivalently, strong topology of $L^2(\Omega;\sph^2)$), we also obtain the convergence of (almost) minimizers as outlined, e.g., in \cite[Proposition 19.11]{BrDe98}.
\end{remark}


\section{Connection to the trace operator}\label{sec:conn_trace}

After proving sufficient conditions for the existence of suitable extensions, we now turn to the analysis of necessary conditions. To that end, we will relax the problem: we will not deal with extensions for functions living on the whole perforated domain but only on one of the constituent cubes. Our goal is to connect this problem to that of the surjectivity of the trace operator for Sobolev functions between manifolds. Since the limit cases \(p=1\) and \(p=\infty\) are always peculiar in the study of trace operators for Sobolev maps, we will restrict ourselves to \(1<p<\infty\).

Throughout this section, \(\p M\) will always denote a \(C^2\)-regular connected, compact, boundaryless hypersurface of \(\R^d\) of dimension \(d-1\) that, by the Jordan separation theorem (cf.\ \cite[Proposition~6.4]{OuRu09}), bounds an open subset of \(\R^d\). We will understand this open subset, denoted by \(M\), together with its boundary \(\p M\), as a submanifold with boundary of \(\R^d\) of dimension \(d\) and class \(C^2\). According to Section~\ref{sec:prelim} and its notation, there exists \(\veps>0\) such that 
\begin{equation}
M_\veps
\coloneqq \{x\in \R^d:\, 0<\dist(x,M)< \veps\}
\end{equation}
is a well-defined one-sided tubular neighborhood (or collar) of \(\p M\).
We fix such a \(\veps\) and also define \(\hat{M}_\veps \coloneqq M\cup M_\veps\).

\begin{remark}\label{rmk:local}
    In the following, we will sometimes consider the case where \(M=B^{d-1}\times [0,1)\) and \(\p M = B^{d-1}\times \{0\} \simeq B^{d-1}\), and define \(M_\veps = B^{d-1}\times(-\veps,0)\). Note that this model deviates from the general assumptions about \(M\) and its boundary. It can be interpreted as a "local" model in the sense that, by drilling a small cylinder inside $\hat{M}_\veps$ defined above, we get something that is diffeomorphic to $B^{d-1}\times(-\veps,1)$.
\end{remark}

We note here that all the following proofs, especially those of Propositions~\ref{prop:gen_nec_cond} and~\ref{prop:gen_nec_cond_smallp}, can be adapted to lower dimensional submanifolds with boundary, since we use only local constructions that can be carried out independently of the dimension of the manifold. In this situation, the annulus \(M_\veps\) will only expand in the normal direction of the boundary \(\p M\). In a further step, one can even generalize the results to arbitrary smooth Riemannian manifolds \(M\) with boundary (compare also \cite[Section~7]{HaLi87}). Then \(M_\veps\) has to be defined by embedding \(M\) isometrically into the Euclidean space, for example with the Nash embedding theorem. Furthermore, all results related to the surjectivity of the trace operator for Sobolev mappings that we will refer to in the sequel, hold in general also for a smooth Riemannian manifold \(M\).

Let the target manifold \(N\) be as in Section~\ref{sec:prelim}, but additionally equipped with a Riemannian structure.
For functions on the boundary \(\p M\) and \(0<s<1\), \(1<p<\infty\), we consider the fractional Sobolev space 
\begin{equation}\label{eq:man_val_frac_sob}
	W^{s,p}(\p M;N)
	\coloneqq \{ f \in W^{s,p}(\p M;\R^l) :\ f(x)\in N \text{ for \ae } x\in \Om\},
\end{equation}
where the definition of the Sobolev space $W^{s,p}(\p M;\R^l)$ can be found, e.g., in \cite[Section~6.7]{KuJoFu77}.

We now define the trace spaces with respect to the target manifold \(N\)
\begin{align*}
\mathcal{T}^p(\p M;N)
&\coloneqq \{f\in W^{1-1/p,p}(\p M;N):\ \exists F\in W^{1,p}(M;N) \text{ s.t. } \Tr(F)= f \text{ on } \p M\},\\
\mathcal{T}^c(\p M;N)
&\coloneqq \{f\in C(\p M;N):\ \exists F\in C(M;N) \text{ s.t. } F = f \text{ on } \p M\},
\end{align*} 
and
\begin{equation*}
\mathcal{E}^p(M_\veps ;N)
\coloneqq \{f\in W^{1,p}(M_\veps;N):\ \exists \hat{F}\in W^{1,p}(\hat{M}_\veps;N) \text{ s.t. } \hat{F}= f \text{ on } M_\veps\}.
\end{equation*} 
Here, \(\Tr(F)\) denotes the image of \(F\in W^{1,p}(M;N)\) in \(W^{1-1/p,p}(\p M;N)\) under the trace operator
\begin{equation}\label{eq:trace_op}
\Tr: W^{1,p}(M;N) \to W^{1-1/p,p}(\p M;N),
\end{equation}
which is the restriction of the usual trace operator \(W^{1,p}(M;\R^l) \to W^{1-1/p,p}(\p M;\R^l)\) to \(W^{1,p}(M;N)\).
We simplify the setting of Section~\ref{sec:constr_ext}, and reformulate our problem as follows:
\renewcommand{\theequation}{\(P_{ex}\)}
\begin{itemize}
	\item[(\(Q_{ex}\))]\label{item:qex}
	Under which assumptions on \(M\), \(N\), and \(p\) does 
	\begin{equation}\label{eq:extension_property}
		\mathcal{E}^p(M_\veps ;N)= W^{1,p}(M_\veps;N)
	\end{equation}
	hold?
\end{itemize}
In other words, we ask in which cases a function that is defined in a small neighborhood around \(M\) and takes values in the manifold \(N\) can be extended to a function on all of \(M\), still satisfying the same manifold constraint. Note that as soon as property \eqref{eq:extension_property} holds for some \(\veps'\le\veps\), the same is automatically true for any parameter in the range \((0,\veps]\), since it is always possible to extend functions in \(W^{1,p}(M_{\veps'};N)\) to \(W^{1,p}(M_\veps;N)\) (for instance, by usual reflection methods).
If, instead, we start with a function given on the boundary \(\p M\), this leads to the following question:
\renewcommand{\theequation}{\(P_{tr}\)}
\begin{itemize}
	\item[(\(Q_{tr}\))]\label{item:qtr}
	Under which assumptions on \(M\), \(N\), and \(p\) does 
	\begin{equation}\label{eq:surjectivity_ptrace}
		\mathcal{T}^p(\p M;N)= W^{1-1/p,p}(\p M;N)
	\end{equation}
	hold?
\end{itemize}
The latter problem corresponds to the surjectivity of the trace operator for manifold-valued Sobolev maps, and it has been extensively addressed in \cite{HaLi87,BeDe95,Be14} and more recently in \cite{MiVS21,MaVS23, VS24}. It turns out that the answer to question (\(Q_{tr}\)) is crucially dependent on the relation between \(p\) and \(d\).

For \(p\ge d\), question (\(Q_{tr}\)) is equivalent (cf.\ \cite[Theorem~1, Theorem~2]{BeDe95}, see also Theorem~\ref{thm:trace_pblm_equiv} below) to finding conditions on \(M\) and \(N\) that ensure the surjectivity of the trace operator for continuous functions, that is, to the following question:
\renewcommand{\theequation}{\(P_{tr}^c\)}
\begin{itemize}
	\item[(\(Q_{tr}^c\))]\label{item:qtr_c}
	Under which assumptions on \(M\) and \(N\) does 
	\begin{equation}\label{eq:surjectivity_ctrace}
		\mathcal{T}^c(\p M;N)= C(\p M;N)
	\end{equation}
	hold?
\end{itemize}
\renewcommand{\theequation}{\arabic{section}.\arabic{equation}}
Question (\(Q_{tr}^c\)) is very well investigated in algebraic topology: a branch usually referred to as {\it obstruction theory} (cf.\ \cite[pp.~415]{Ha02} and \cite[pp.~175]{Hu59}) addresses the question of extendability of continuous functions in terms of the topologies of \(M\) and \(N\).

In the case \(p<d\), the answer to question (\(Q_{tr}\)) depends in an intricate way on \(p\) and on the topology of \(N\).
It is now well-known that, if \(N\) is \(\lfloor p-1\rfloor\)-connected, the trace operator is surjective (see \cite[Theorem~6.2]{HaLi87}, whose proof relies on the retraction from Lemma~\ref{lem:lip_retract}).
In recent years it has been shown that several obstructions can appear that prevent the trace operator from being surjective (see, e.g., \cite{MaVS23}).
However, in the special case that \(1\le p<4\) and for the special choice \(M=B^{d-1}\times [0,1)\) with \(\p M = B^{d-1}\times \{0\}\) (see Remark~\ref{rmk:local}), question (\(Q_{tr}\)) was completely characterized in terms of the homotopy groups of the target manifold (see \cite[Corollary~1.1]{MiVS21}). 
The recent preprint \cite{VS24} provides an answer to question (\(Q_{tr}\)) for all \(p<d\) and all manifolds \(M\). It states that a necessary and sufficient condition for property \eqref{eq:surjectivity_ptrace} to hold is that the homotopy groups \(\pi_1(N),\dots, \pi_{\lfloor p-2\rfloor}(N)\) are finite and an additional topological condition on \(M\) and \(N\) is fulfilled to prevent the appearance of topological obstructions.

In the remaining part of this section we elaborate on how properties \eqref{eq:extension_property} and \eqref{eq:surjectivity_ptrace} are related to one another. First note that, since every function in \(W^{1,p}(M_\veps;N)\) possesses a trace in \(W^{1-1/p,p}(\p M;N)\), the following implication is immediate.

\begin{proposition}\label{prop:suff_cond}
    Let \(1< p< \infty\). As soon as property \eqref{eq:surjectivity_ptrace} holds, then also property \eqref{eq:extension_property} is fulfilled.
\end{proposition}

In contrast, the necessity of property \eqref{eq:surjectivity_ptrace} for property \eqref{eq:extension_property} is not as obvious. In general, one can ask:

\begin{center}
\parbox{.9\textwidth}{Is there additional information that allows to extend any function given on the annulus \(M_\veps\) to the inside of \(M\), even though the trace operator is not surjective?}
\end{center}

For example, when \(\p M = B^{d-1}\times \{0\}\), \(M=B^{d-1}\times [0,1)\), and \(M_\veps=B^{d-1}\times(-\veps,0)\) (see Remark~\ref{rmk:local}), every function on \(M_\veps\) can be extended to \(\hat{M}_\veps\) by reflection (independently of the topological properties of $N$), so that property \eqref{eq:extension_property} holds; but, at the same time, in the case $2 \le p<d$, if \(\pi_{\lfloor p-1\rfloor}(N)\ne 0\) then property \eqref{eq:surjectivity_ptrace} turns out to be false (see \cite[Theorem~3]{MiVS21}). Note that, in such a case, those functions in the space $W^{1-1/p,p}(\p M;N)$ that cannot be extended to \(M\) are also not traces of functions in \(W^{1,p}(M_\veps;N)\) (otherwise the trace operator would be surjective).

\begin{remark}
    The above question is specific to Sobolev maps between two manifolds. For Sobolev maps into the Euclidean space, properties \eqref{eq:surjectivity_ptrace} and \eqref{eq:extension_property} are always equivalent because the problems of extending a function to the inside or the outside of \(\p M\) cannot be distinguished due the absence of the target manifold condition.
\end{remark}

We have the following easy lemma, that gives a sufficient condition for the necessary condition given by property \eqref{eq:surjectivity_ptrace}:

\begin{lemma}\label{lem:SCfortheNC}
    Assume that property \eqref{eq:extension_property} holds. As soon as it is possible to extend a given function in \(W^{1-1/p,p}(\p M;N)\) to a function in \(W^{1,p}(M_\veps;N)\), we can infer that also property \eqref{eq:surjectivity_ptrace} is fulfilled.
\end{lemma}

In the case \(p\ge d\), it turns out that the obstructions for question (\(Q_{tr}\)), which could also depend on the topology of \(M\), do not prevent us from expanding fractional Sobolev functions from the boundary to the exterior annulus. Note that, since \(M_\veps\) does not posses any additional hole, \(\p M\) and \(M_\veps\) share the same topological properties. This feature is implicitly used in Proposition~\ref{prop:gen_nec_cond} below, proving that, indeed, in this case we can establish the necessity of property \eqref{eq:surjectivity_ptrace} for property \eqref{eq:extension_property}, and that both problems are equivalent. 

In contrast, for \(p<d\), the problem of extending functions in \(W^{1-1/p,p}(\p M;N)\) to \(W^{1,p}(M_\veps;N)\) is much more sensitive to the topology of \(N\). However, it is clear that as soon as \eqref{eq:surjectivity_ptrace} is independent of the topology of \(M\) and the target manifold allows to extend functions in \(W^{1-1/p,p}(\p M;N)\) to the inside of \(M\), they can also be extended to the outside. Correspondingly, as soon as one of the above is not possible, nor is the other. In these instances, the answer to question (\(Q_{tr}\))
is closely related to the necessity of property \eqref{eq:surjectivity_ptrace} for property \eqref{eq:extension_property}. 
For example, taking into account Proposition~\ref{prop:suff_cond}, as soon as the conditions on \(M\), \(N\), and \(p\) are met such that property \eqref{eq:surjectivity_ptrace} 
is fulfilled, then properties \eqref{eq:surjectivity_ptrace} and \eqref{eq:extension_property} naturally coincide -- and hold simultaneously. We validate this observation below in Proposition~\ref{prop:gen_nec_cond_smallp} for the case where \(N\) is \(\lfloor p-1\rfloor\)-connected.
In Proposition~\ref{prop:case2.2} we extend this result to target manifolds $N$ for which \(\pi_1(N)\) is only finite, by additionally assuming that \(\p M\) is simply connected.
Again, this will turn out to be a case where properties \eqref{eq:surjectivity_ptrace} and \eqref{eq:extension_property} are valid.
 
To our knowledge, in all the remaining cases the question of whether property \eqref{eq:surjectivity_ptrace} is necessary for property \eqref{eq:extension_property} remains completely open at the moment. 
 
Let us now make our comments rigorous. The case distinction implicitly applied above is inspired by the analysis of question (\(Q_{tr}\)) in \cite{Be14} and \cite{MiVS21}.

\subsection{The case \texorpdfstring{\(p\ge d\)}{p>d}}\label{subsec:pged}
We begin with a short standard example. 

Take \(M\coloneqq\ol{B}^d\) the \(d\)-dimensional closed unit ball with boundary \(\p M=\Sp^{d-1}\) and \(N\coloneqq\Sp^{d-1}\). We further define \(M_\veps=\{x\in \R^d:1<|x|<1+\veps\}\) the annulus around \(M\). Then the function 
\begin{equation}\label{eq:counterex}
f: M_\veps\to N,\quad
x\mapsto \frac{x}{|x|}
\end{equation} 
is in \(W^{1,p}(M_\veps;N)\) for any \(p\in [1,\infty]\). 
A simple topological argument will show that, at least for \(p\ge d\), the function $f$ cannot posses an extension in \(W^{1,p}(M;N)\).
In \cite[Theorems~1 and 2]{BeDe95} the following result is proved. 
\begin{theorem}\label{thm:trace_pblm_equiv}
    Let \(p\ge d\). Then property \eqref{eq:surjectivity_ptrace} is true if and only if property \eqref{eq:surjectivity_ctrace} holds.	 	
\end{theorem}
For the special choice \(M=\ol{B}^d\), it is well-known that property \eqref{eq:surjectivity_ctrace} holds if and only if \(\pi_{d-1}(N)=\{0\}\) (see, e.g., \cite[p.~346]{Ha02}). Combining this observation with the above Theorem~\ref{thm:trace_pblm_equiv}, we obtain the following:
\begin{corollary}\label{cor:trace_pblm_equiv2}
	Let \(p\ge d\). Then property \eqref{eq:surjectivity_ptrace} is true for \(M=\ol{B}^d\) if and only if
	\begin{equation}\label{eq:ext_con}
	\pi _{d-1}(N)=\{0\}.
	\end{equation}
\end{corollary}
Let us focus on the case $N = \Sp^{d-1}$. It is known that \(\pi_{d-1}(\Sp^{d-1})=\Z\). From Corollary~\ref{cor:trace_pblm_equiv2} it follows that property \eqref{eq:surjectivity_ptrace} does not hold, that is, there exists (at least) a function in the space $W^{1-1/p,p}(\Sp^{d-1};\Sp^{d-1})$ that cannot be extended to the interior of $\ol{B}^d$. We are going to show that an example is provided exactly by the function defined in \eqref{eq:counterex}.
Indeed, since $f$ is such that \(f|_{\Sp^{d-1}}=\id\), it has Brouwer-Kronecker degree \(\deg(f|_{\Sp^{d-1}})=1\). On the other hand, any continuous extension to \(\ol{B}^d\) should have degree equal to zero (cf.\ \cite[Propostion~4.4]{OuRu09}). We conclude that \(f\) cannot admit an extension in \(W^{1,p}(\ol{B}^d;\Sp^{d-1})\) for \(p\ge d\), since any such extension would necessarily have trace equal to the identity on \(\p M\). 

The idea underlying the previous example is that the non-surjectivity of
the trace operator in \eqref{eq:trace_op} can pose an obstruction to the extension of functions in \(W^{1,p}(M_\veps;N)\) to the inside of \(M\). We prove below the following general result in that direction.
\begin{proposition}\label{prop:gen_nec_cond}
	Let \(p\ge d\). If property \eqref{eq:extension_property} holds, then also property \eqref{eq:surjectivity_ptrace} is fulfilled and the trace operator in \eqref{eq:trace_op} is surjective.
\end{proposition}

As a corollary, from Proposition~\ref{prop:gen_nec_cond} together with Theorem~\ref{thm:trace_pblm_equiv}, we get a necessary condition in terms of the surjectivity of the trace operator
\begin{equation}\label{eq:trace_op_cont}
\Tr: C(M;N) \to C(\p M;N)
\end{equation}
for continuous functions.
\begin{corollary}
	Let \(p\ge d\). If property \eqref{eq:extension_property} holds, then also property \eqref{eq:surjectivity_ctrace} is fulfilled and the trace operator in \eqref{eq:trace_op_cont} is surjective.
\end{corollary}

Together with Proposition~\ref{prop:suff_cond} and again Theorem~\ref{thm:trace_pblm_equiv}, which establish the sufficiency of property \eqref{eq:surjectivity_ctrace} for property \eqref{eq:extension_property}, this shows that for \(p\ge d\) all properties \eqref{eq:extension_property}, \eqref{eq:surjectivity_ptrace}, and \eqref{eq:surjectivity_ctrace} are actually equivalent.

\begin{proof}[Proof of Proposition~\ref{prop:gen_nec_cond}]
Let \(p\ge d\). We need to prove that if every function in \(W^{1,p}(M_{\veps};N)\) can be extended to a function in \(W^{1,p}(\hat{M}_\veps;N)\), then also every function in the space \(W^{1-1/p,p}(\p M;N)\) can be extended to \(M\) to a function in \(W^{1,p}(M;N)\).

Let a function \(f\in W^{1-1/p,p}(\p M;N)\) be given. From Lemma~\ref{lem:SCfortheNC}, it suffices to prove that there exists an extension \(\hat{F}\in W^{1,p}(M_\veps;N)\) with trace \(f\) on \(\p M\) (at least for \(\veps'\le \veps\) sufficiently small, compare \eqref{eq:choice_veps}). Then, extending \(\hat{F}\) to \(M\) will also yield the desired extension for \(f\). The construction is based on standard arguments that are usually used to directly prove the surjectivity of the trace operator in \eqref{eq:trace_op} (see, for instance, \cite{Leoni17, KuJoFu77}).

\noindent
\textit{Step 1.}
Let \((U_i,\varphi_i)_{i\in I}\), where \(I=\{1,\cdots,\mathfrak{m}\}\), be a finite atlas of the compact \(C^2\) manifold \(\p M\) with, for all \(i\in I\), open subsets \(U_i\subset \p M\) and local charts \(\varphi_i:U_i\to \varphi(U_i)\subset \R^{d-1}\) that are diffeomorphisms. We set
\begin{equation}
	V_i\coloneqq \varphi_i(U_i) 
	\quad \text{and} \quad
	\phi_i\coloneqq\varphi_i^{-1}: V_i\to U_i.
\end{equation}
Let further \((\zeta_i)_{i\in I}\) be a partition of unity associated to the covering \((U_i)_{i\in I}\) of \(\p M\). 
By \(\nu:\p M\to \Sp^{d-1}\) we denote the continuous normal vector field associated to \(\p M\). We define the maps
\begin{equation}
	\psi_i: V_i\times [0,\veps]\to M_\veps,
	\quad (x',h) \mapsto \psi_i(x',h) \coloneqq \phi_i(x')+h\nu(\phi_i(x'))
	\quad \forall i\in I.
\end{equation}
Then all \(\psi_i\) are diffeomorphisms with \(\psi_i(V_i\times \{0\})=U_i\subset \p M\). We define the pulled back functions \(g_i:V_i\to N\) as the composition in local coordinates given by
\begin{equation}\label{eq:gi}
	g_i(x')\coloneqq(\zeta_if)(\psi_i(x',0)) 
	\quad \text{for } x'\in V_i.
\end{equation}
It follows readily that \(g_i\in W^{1-1/p,p}(V_i;N)\).
We denote by \(B^{d-1}_1\) the unit ball in \(\R^{d-1}\) and by \(\eta\in C_c^\infty(B^{d-1}_1)\) a standard mollifier with \(\norm{\eta}{L^1(\R^{d-1})}=1\), and set \(\eta_h(x)\coloneqq h^{-(d-1)}\eta(x/h)\). Without loss of generality, we can take \(\eta\) to be positive. We might now extend the local functions \(g_i\) through the convolution 
\begin{equation}
	g_i^{ex}(x'')
	\coloneqq \int_{B^{d-1}_h(x')} \eta_h\left(x'-y'\right)g_i(y') {\dl y'}, \quad x'' \coloneqq (x',h) \in V_i \times \R_+.
\end{equation}
It can be proved that \(g_i^{ex}\in C^1(V_i\times\R_+;\R^l)\) for all \(i\in I\) (cf.\ \cite[Theorem~18.28]{Leoni17}). It follows that
\begin{equation}\label{eq:sob_extension}
	f^{ex}(\xi + h\nu(\xi))
	\coloneqq \sum_{i\in I} \left( g_i^{ex}\circ\psi_i^{-1}\right) (\xi + h\nu(\xi))
	\quad \text{for } \xi \in \p M, \ h \in [0,\veps],
\end{equation}
defines a global extension of \(f\) to a function in \(W^{1,p}(M_\veps; \R^l)\). For more details we refer to \cite[Theorem~18.40]{Leoni17}.

\noindent
\textit{Step 2.}
It remains to prove that, even though for the critical case \(p=d\) neither \(f\) on \(\p M\) nor \(f^{ex}\) in \(M_\veps\) are necessarily continuous, the extension \(f^{ex}\) will always take values in \(N_\delta\) if we restrict it to a one-sided tubular neighborhood of \(\p M\) of small enough thickness. Note that a similar strategy is applied in \cite[Theorem~2]{BeDe95} to prove the equivalence between properties \eqref{eq:surjectivity_ptrace} and \eqref{eq:surjectivity_ctrace}. For convenience, we repeat the calculations at this point in some detail. We set
\begin{equation*}\label{eq:convention}
	f(\psi_i(\,\cdot\,,0)) = 0 
	\quad \text{ in } \R^{d-1}\setminus V_i.
\end{equation*}
For \(\xi \in \p M\), let $x_i'=\varphi_i(\xi)$ and $x_i''= (x_i',h) = (\varphi_i(\xi),h)$.
Let also $\xi^0 \in \p M$ and $z_i' = \varphi_i(\xi^0)$. Then, in light of the properties of the mollifier \(\eta\) and the partition of unity \((\zeta_i)_{i\in I}\), from \eqref{eq:sob_extension} and \eqref{eq:gi} we infer
\begin{align*}
    \begin{split}
        f^{ex}(\xi &+ h\nu(\xi)) - f(\xi^0) = \sum_{i\in I} \Big(\big( g_i^{ex}\circ\psi_i^{-1}\big) (\xi + h\nu(\xi)) - \zeta_i(\xi)f(\xi^0)\Big) \\
		&= \sum_{i\in I} \Bigg( g_i^{ex}(x_i'') - \zeta_i(\phi_i(x_i'))\sum_{j\in I}(\zeta_j f)(\phi_j(z_j'))\Bigg).
    \end{split}
\end{align*}
Note that the product $\zeta_i(\phi_i(x_i'))\zeta_j(\phi_j(z_j'))$ vanishes whenever \(\xi\in U_i\) and \(\xi^0\in U_j\setminus(U_j \cap U_i)\) for \(i\neq j\). Moreover, $\xi$ and $\xi^0$ can both lie in two (or more) distinct $U_i$'s. We therefore estimate
\begin{align}\label{eq:difference}
    \begin{split}
        &f^{ex}(\xi + h\nu(\xi)) - f(\xi^0) \\
		&\le \sum_{i\in I} \int_{B^{d-1}_h(x'_i)} \eta_h(x'_i-y')\Big(g_i(y')-\zeta_i(\phi_i(x_i'))\zeta_i(\phi_i(z_i'))f(\phi_i(z_i'))\Big) {\dl y'} \\
		&= \sum_{i\in I}\Bigg( \int_{B^{d-1}_h(x'_i)} \eta_h(x'_i-y')\zeta_i(\phi_i(x_i'))\zeta_i(\phi_i(z_i'))
		\Big( f(\phi_i(y'))-f(\phi_i(z_i'))\Big) {\dl y'}
		+ R_i(z_i',h) \Bigg),
	\end{split}
\end{align}
where the remainders are given by
\begin{equation}\label{eq:remainder}
R_i(z_i',h) 
= \int_{B^{d-1}_h(x'_i)} \eta_h(x'_i-y')\Big(\zeta_i(\phi_i(y')-\zeta_i(\phi_i(x_i'))\zeta_i(\phi_i(z_i'))\Big)f(\phi_i(y')) {\dl y'}.
\end{equation}
We postpone the asymptotic analysis of the remainder and focus on the first contribution on the right-hand side of \eqref{eq:difference}.
We may now assume, for the sake of simplicity, that \(U_i\) is already flat and contained in \(\R^{d-1}\times \{0\}\). Thus, we can take the chart \(\varphi_i\) to be the identity. For \(x'_i\in U_i\subset \R^{d-1}\times \{0\}\) and $z' \in B^{d-1}_h(x_i')$ we define
\begin{align*}
    \begin{split}
		F_i(z',h) &\coloneqq \int_{B^{d-1}_h(x'_i)} \eta_h(x'_i-y')\zeta_i(\phi_i(x_i'))\zeta_i(\phi_i(z'))
		\Big( f(\phi_i(y'))-f(\phi_i(z'))\Big) {\dl y'} \\
		&= \int_{B^{d-1}_h(x'_i)} \eta_h\left(x'_i-y'\right)\zeta_i(x_i')\zeta_i(z')\Big(f(y')-f(z')\Big) {\dl y'}.
    \end{split}
\end{align*}
Notice that
\begin{equation*}\label{eq:simple_est}
|y'-z'|
\le 2h
\quad \text{for all } y' \in B^{d-1}_h.
\end{equation*}
Hence, by Jensen's inequality and the fact that \(0\le \zeta_i\le 1\) for all $i \in I$, we get the estimate
\begin{align}\label{eq:slobodeckij}
	\begin{split}
	&\frac{1}{\mathcal{L}^{d-1}(B^{d-1}_h(x'_i))}\int_{B^{d-1}_h(x'_i)} |F_i(z',h)|^p {\dl z'}\\
	&\le C\frac{\norm{\eta}{L^\infty(\R^{d-1})}^p}{h^{d-1}\mathcal{L}^{d-1}(B^{d-1}_h(x'_i))}\Bigg(
	\int_{B^{d-1}_h(x'_i)}\int_{B^{d-1}_h(x'_i)} |f(y')-f(z')|^p {\dl y'}{\dl z'} \Bigg) \\
	&\le C\norm{\eta}{L^\infty(\R^{d-1})}^p\Bigg(h^{p-d}
	\int_{B^{d-1}_h(x'_i)}\int_{B^{d-1}_h(x'_i)} \frac{|f(y')-f(z')|^p}{|y'-z'|^{d+p-2}} {\dl y'}{\dl z'} \Bigg)
	\end{split}
\end{align}
for some constant \(C>0\) independent of \(h\). Since the double integral on the right-hand side coincides with the Gagliardo semi-norm of \(f\) on \(B_h^{d-1}(x'_i)\) for \(s=1-1/p\) (see \cite{Ga57}), it decays to zero for \(h\to 0\) by the absolute continuity of the Lebesgue integral. Consequently, by the mean value theorem, and since we are assuming $p \ge d$, there exist some threshold \(\veps_1>0\) and some \(z'\in B^{d-1}_h(x'_i)\) such that, for all $i \in I$,
\begin{equation*}
	|F_i(z',h)|^p
	\le \left(\frac{\delta}{2\mathfrak{m}}\right)^p
	\qquad \forall h\le \veps_1.
\end{equation*}
Since \((\varphi_i)_{i\in I}\) is a finite family of diffeomorphisms, the same estimate holds true, up to possibly altering \(\veps_1\), also for general \(U_i\subset \p M\), uniformly in \(i\in I\) and \(\xi\in \p M\).
As for the remainder in \eqref{eq:remainder}, we recall that \(f\in L^\infty(\p M;N)\) because \(N\) is compact. A straightforward computation shows that
\begin{align*}
	|R_i(&z_i',h)|
	\le \int_{B^{d-1}_h(x'_i)} \eta_h(x'_i-y')\big|\zeta_i(\phi_i(y')-\zeta_i(\phi_i(x_i'))\zeta_i(\phi_i(z_i'))\big|\big|f(\phi_i(y'))\big| {\dl y'}\\
	&\le \mathcal{L}^{d-1}(B^{d-1})\norm{\eta}{L^\infty(\R^{d-1})}\norm{f}{L^\infty(\p M)}\norm{\zeta_i\circ\phi_i-\zeta_i(\phi_i(x_i'))\zeta_i(\phi_i(z_i'))}{L^\infty(B^{d-1}_h(x'_i))}.
\end{align*}
The expression \(\norm{\zeta_i\circ\phi_i-\zeta_i(\phi_i(x_i'))\zeta_i(\phi_i(z_i'))}{L^\infty(B^{d-1}_h(x'_i))}\) tends to zero for \(h\to 0\) uniformly in \(i\in I\) and \(z_i'\), since \((\zeta_i)_{i\in I}\) is a finite family of smooth partition functions. Thus, there exists \(\veps_2>0\) such that
\begin{equation}
	|R_i(z_i',h)|
	\le \frac{\delta}{2\mathfrak{m}}
	\qquad \forall h\le \veps_2.
\end{equation}
Overall, we obtain the existence of some \(\veps'>0\) and, for almost every \(\xi \in \p M\), of some \(\xi^0\in \p M\) with \(f(\xi^0)\in N\) such that
\begin{equation}\label{eq:choice_veps}
	\dist(f^{ex}(\xi + h\nu(\xi)),N)
	\le |f^{ex}(\xi + h\nu(\xi)) - f(\xi^0)|
	\le \delta
	\qquad\forall h\le \veps'.
\end{equation}
Consequently, we find that \(f^{ex}(x)\in N_\delta\) for \ae \(x\in M_{\veps'}\). Without loss of generality, we may assume that \(\veps'= \veps\). Otherwise, if \(\veps'<\veps\), we can now extend \(f^{ex}\) to the rest of \(M_\veps\) (e.g., by suitable reflections). Hence, we have that \(f^{ex}\in W^{1,p}(M_{\veps};N_\delta)\).

\noindent
\textit{Step 3.}
For \(\xi\in C^1(\ol{N_\delta};N)\) the nearest point projection of Lemma~\ref{lem:nearest_pp}, we have immediately that \(\xi\circ f^{ex}\in W^{1,p}(M_{\veps};N)\). By the assumption that \eqref{eq:extension_property} is true, we are able to find some \(\hat{F}\in W^{1,p}(\hat{M}_\veps;N)\) that coincides with \(\xi\circ f^{ex}\) on the annulus \(M_\veps\). We finally define 
\begin{equation}
	F\coloneqq\hat{F}|_{M}
\end{equation}
as the candidate for the extension of \(f:\p M\to N\) to the inside of \(M\). This proves the claim, since by construction there holds indeed \(F|_{\p M}=f\).
\end{proof}

\begin{remark}
    An alternative (but less self-contained) proof of Proposition~\ref{prop:gen_nec_cond} can be obtained by using the continuous embedding of \(W^{1-1/p,p}(\p M;N)\) into the space \(\VMO(\p M;N)\) of functions of vanishing mean oscillation (compare for instance \cite[Section~1.2]{BrNi95}) and noting that, for \(f^{ex}\) as in \eqref{eq:sob_extension}, inequality \eqref{eq:choice_veps} follows from a similar property of \(\VMO\) functions (cf.\ inequality (7) on p.~206 in \cite{BrNi95}).
\end{remark}

\subsection{The case \texorpdfstring{\(p<d\)}{p<d} and \texorpdfstring{\(N\)}{N} \texorpdfstring{\(\lfloor p-1\rfloor\)}{[p-1]}-connected}
Unfortunately, the same strategy adopted for Proposition~\ref{prop:gen_nec_cond} will not be generally applicable for \(p<d\). One exception is given for example when the target manifold is \(\lfloor p-1\rfloor\)-connected. In this case, a result similar to Proposition~\ref{prop:gen_nec_cond} can be deduced by replacing the nearest point projection with the generalized retraction from Lemma~\ref{lem:lip_retract}, and proceeding as in the proof of Theorem~\ref{thm:target_ext}. This result can be seen as the counterpart of \cite[Theorem~6.2]{HaLi87}, where the authors proved that under the same \(\lfloor p-1\rfloor\)-connectedness assumption any fractional Sobolev function on $\p M$ can be extended to the inside of \(M\).
In other words, under the special \(\lfloor p-1\rfloor\)-connectedness condition, there always exists a retraction onto \(N\), and thus, not only properties \eqref{eq:surjectivity_ptrace} and \eqref{eq:extension_property} become again equivalent, but as in the unconstrained case, where no target manifold is prescribed, both problems are also solvable simultaneously.

\begin{proposition}\label{prop:gen_nec_cond_smallp}
	Let \(p< d\). Let \(N\) be as in Section~\ref{sec:prelim} and additionally \(\lfloor p-1\rfloor\)-connected. If property \eqref{eq:extension_property} holds, then also property \eqref{eq:surjectivity_ptrace} is fulfilled and the trace operator in \eqref{eq:trace_op} is surjective. Moreover, in this case properties \eqref{eq:surjectivity_ptrace} and \eqref{eq:extension_property} are simultaneously fulfilled.
\end{proposition}

\begin{proof}
    Let \(f\in W^{1-1/p,p}(\p M;N)\), and let \(f^{ex}\in W^{1,p}(M_{\veps}; \R^l)\) be its extension to $M_\veps$ given by formula \eqref{eq:sob_extension}. Note that, by construction, \(f^{ex}\) will only take values in some compact subset of \(\R^l\). We are now going to adapt the steps from the proof of Theorem~\ref{thm:target_ext}.
    This time we approximate \(f^{ex}\) by a sequence \((f_k)_k\subset C^\infty(\ol{M_{\veps}};\R^l)\) of smooth functions converging strongly to $f^{ex}$ in $W^{1,p}(M_\veps;\R^l)$, and such that each of these functions again ranges in a compact subset of \(\R^l\). Note that, due to the regularity of \(\p M\) (and then of $\p M_\veps$), we can actually take these approximations to be smooth up to the boundary (cf., for instance, \cite[Theorem~11.35]{Leoni17}, or also \cite[Section~5.3.3, Theorem~3]{Ev98}). This entails that, as in the proof of Theorem~\ref{thm:target_ext}, we can find a function \(\tilde{f}\in W^{1,p}(M_{\veps};N)\) such that
	\begin{equation}
		(P_{h_k}|_N)^{-1}\circ P_{h_k} \circ f_k
		\rightharpoonup \tilde{f} \ \mbox{ in } W^{1,p}(M_{\veps};\R^l),
	\end{equation}
	where the \(P_{h_k}\) are the translated retractions defined in \eqref{eq:translated_proj}. It remains to prove that 
	\begin{equation}\label{eq:trace_prop}
		\Tr(\tilde{f}) 
		= f,
	\end{equation} 
	where \(\Tr: W^{1,p}(M;N)\to W^{1-1/p}(\p M;N)\) is the trace operator from \eqref{eq:trace_op}. Since the trace operator is defined as the continuous extension of the operator restricting functions in \(C_c^\infty(\R^d)\) to \(\p M\), we know that \(f= \Tr(f^{ex})\) is the strong limit of \((f_k|_{\p M})_k\) in \(L^p(\p M)\).
	It follows that, up to some subsequence that we do not relabel, the maps \(f_k\) converge pointwise \ae (with respect to the (\(d-1\))-dimensional Hausdorff measure \(\mathcal{H}^{d-1}\)) to \(f\) on \(\p M\). Besides, the local Lipschitzianity of the retractions \(P_{h_k}:(Q_R+h_k)\setminus( X+h_k)\to N\) (compare Lemma~\ref{lem:lip_retract}) guarantees that
	\begin{equation}
		(P_{h_k}|_N)^{-1}\circ P_{h_k} \circ f_k 
		\in C(\ol{M_{\veps}}\setminus f_k^{-1}(X+h_k);N)
		\quad \forall k\in \N.
	\end{equation}
In particular, since \((N_\delta+h_k)\cap (X+h_k)=\emptyset\) (see again Lemma~\ref{lem:lip_retract}) and by construction we have that \(h_k<\delta\) for \(k\) large enough, the composition is well-defined for \(x\in \ol{M_{\veps}}\) as soon as \(f_k(x)\) is sufficiently close to \(N\subset N_\delta+h_k\). Hence, we find that by choosing an approximating sequence \(\sequence{f}{k}\) of functions that are smooth up to the boundary \(\p M\), the pointwise characterization of the limit function in \eqref{eq:def_ext1} of Corollary~\ref{cor:aux_comp} is actually valid up to the boundary. To be precise, by the pointwise convergence of the sequence \((f_k)_k\) to \(f\) on \(\p M\) we find that 
\begin{equation}
    \tilde{f}(x)
    =\lim_{k \to \infty} \Big((P_{h_k}|_N)^{-1} \circ P_{h_k} \circ f_k(x)\Big) 
    = (P_{\tilde{h}}|_N)^{-1}\circ P_{\tilde{h}} \circ f(x)
    = f(x)
\end{equation}
for \ae \(x\in \p M\), where \(\tilde{h}<\delta\) is the limit of (a subsequence of) the sequence \((h_k)_k\). This implies that \(\tilde{f}\) coincides \ae on \(\p M\) with \(f\), so that \eqref{eq:trace_prop} is proved. As in step 3 of the proof of Proposition~\ref{prop:gen_nec_cond}, since property \eqref{eq:extension_property} is fulfilled by assumption, there exists now \(\hat{F}\in W^{1,p}(\hat{M}_\veps;N)\) such that \(\hat{F}=\tilde{f}\) on the annulus \(M_\veps\). The function \(F\coloneqq\hat{F}|_{M}\in W^{1,p}(M;N)\) is then the claimed extension of \(f\) to \(M\).

Finally, note that the above strategy also provides an extension \(F\in W^{1,p}(M;N)\) of \(f\) such that $\Tr(F)= f$ on $\p M$ (see also \cite[Theorem~6.2]{HaLi87}). The last part of the claim is then a direct consequence of Proposition~\ref{prop:suff_cond}.
\end{proof}

\subsection{The case \texorpdfstring{\(p<d\)}{p<d} and \texorpdfstring{\(N\)}{N} not simply connected}
Following the latter Proposition~\ref{prop:gen_nec_cond_smallp}, for the last part of this section we will assume that \(2\le p<d\); otherwise, \(\pi_1(N)\) and higher homotopy groups have no relevance for the extension of functions from \(\p M\) to \(M_\veps\). If \(2\le p< 3\le d\), then property \eqref{eq:surjectivity_ptrace} is not fulfilled if \(\pi_1(N)\) is not trivial, as shown in \cite[Theorem~4]{BeDe95} (see also \cite[Section~6.3]{HaLi87} and \cite[Theorem~3]{MiVS21}). Since the same proof can be carried out with \(M_\veps\) instead of \(M\), we cannot deduce any information about whether property \eqref{eq:extension_property} is fulfilled (see Lemma~\ref{lem:SCfortheNC}). This narrows down the range of the integrability exponent to \(3\le p<d\). 

Moreover, if \(\pi_1(N)\) is infinite, it follows plainly from \cite[Theorem~5]{MiVS21} (see also Section~5 in \cite{Be14}) or \cite[Lemma~1.3]{Be14} (with \(M\) replaced by \(M_\veps\) and for $\p M$ simply connected) that we cannot hope to extend functions from \(\p M\) to \(M_\veps\). The same negative result actually holds if any of the homotopy groups \(\pi_1(N), \dots, \pi_{\lfloor p-1\rfloor}(N)\) is infinite.

On the positive side, we have the following result.
\begin{proposition}\label{prop:case2.2}
    Let \(3\le p<d\) and the manifold $M$ be such that \(\p M\) is simply connected. Then, if \(\pi_1(N)\) is finite and 
    \begin{equation*}
        \pi_2(N)=\cdots=\pi_{\lfloor p-1\rfloor}(N)=0,
    \end{equation*}
    property \eqref{eq:extension_property} can only hold if property \eqref{eq:surjectivity_ptrace} also holds. Actually, under the above assumptions, both properties \eqref{eq:surjectivity_ptrace} and \eqref{eq:extension_property} hold.
\end{proposition}
\begin{proof}
    The first part of the statement follows exactly as in the proof of \cite[Theorem~6]{MiVS21}, where it is shown that under the above assumptions, property \eqref{eq:surjectivity_ptrace} is fulfilled.
    Let \(\pi:\tilde{N}\to N\) be the universal covering of \(N\). Notice that \(\tilde{N}\) is compact because \(\pi_1(N)\) is finite. Additionally, the universal covering is simply connected, so that \(\pi_1(\tilde{N})=0\). It can also be shown that \(\pi_j(\tilde{N})=\pi_j(N)=0\) for all \(j=2,\dots,\lfloor p-1\rfloor\) (see, e.g., \cite[Proposition~4.1]{Ha02}). Since \(\tilde{N}\) is compact and \(\p M\) is simply connected, we can apply \cite[Theorem~1]{MiVS21_1} to find for every \(f\in W^{1-1/p,p}(\p M;N)\) a lifting \(\tilde{f}\in W^{1-1/p,p}(\p M;\tilde{N})\) such that \(f=\pi\circ \tilde{f}\). Applying now Proposition~\ref{prop:gen_nec_cond_smallp} with \(N\) replaced by \(\tilde{N}\), we deduce the existence of an extension \(\tilde{F}\in W^{1,p}(\hat{M}_\veps;\tilde{N})\) that coincides with \(\tilde{f}\) on the annulus \(M_\veps\). Set \(\hat{F}\coloneqq\pi \circ\tilde{F}\). Then \(F\coloneqq \hat{F}|_M \in W^{1,p}(M;N)\) is an extension of \(f\).

    The second statement follows again from Proposition~\ref{prop:suff_cond}, this time together with \cite[Theorem~6]{MiVS21} that provides for every function \(f\in W^{1-1/p,p}(\p M;N)\) an extension \(F\in W^{1,p}(M;N)\).
\end{proof}

\section{Conclusions and possible generalizations}\label{sec:concl}

We constructed in Section~\ref{sec:constr_ext} suitable extensions for Sobolev maps defined on a perforated domain and taking values in a compact and $\lfloor p-1\rfloor$-connected manifold $N$ without boundary. The extensions additionally admit bounds independent of the scale, and this allowed us to show in Section~\ref{sec:ex} how to fruitfully employ it in homogenization problems in micromagnetics, where the target manifold is \(N=\sph^2\). The key observation in view of Theorem~\ref{thm:target_ext} is that for dimension \(d=3\) and integrability exponent \(p=2\) there holds \(1\le p<d\) and \(\pi_{\lfloor p-1\rfloor}(N)=\pi_1(\sph^2)=0\).

As examples of target manifold conditions appearing in the literature, we mentioned also $\SL(3)$ and $\SO(3)$ in the introduction. However, in addition to \(\SL(3)\) being not compact, both $\SL(3)$ and $\SO(3)$ are not simply connected. This motivates us to investigate possible enhancements to our construction of extensions. While the compactness issue could be solved by looking at a compactified version of the manifold in the spirit of Alexandroff compactification, the lack of connectedness poses the biggest hurdle. In particular, motivated by the discussion in Section~\ref{sec:conn_trace} about the connection between the extendability of Sobolev functions between manifolds and the surjectivity of the trace operator for such functions, two further research directions seem relevant. First, what can we say in the case that \(2\le p<3\le d\) and \(\pi_1(N)\ne 0\)? We know that the trace operator is not surjective, but this does not necessarily imply that also no extension can exist (see also the discussion before Proposition~\ref{prop:case2.2}). A solution to this problem would include the cases \(N=\SL(3)\) and \(N=\SO(3)\) for \(p=2\) in dimension \(d=3\). However, if one allows for a larger dimension \(d\), keeping the constraint that \(2\le p < d\), then the higher homotopy groups will start to play a role. Thus, inspired by the positive result in Proposition~\ref{prop:case2.2}, the second line of investigation could be whether it is possible to construct an extension if \(\pi_1(N)\) is finite and \(\pi_2(N)=\dots=\pi_{\lfloor p-1\rfloor}(N)=0\). Possible applications could be again \(N=\SL(3)\) or \(N=\SO(3)\) for \(p=3<d\), because \(\pi_1(N)=\Z_2\) and \(\pi_2(N)=0\) in both cases. 

However, despite these special examples, the somehow smoldering question seems to be the one about necessary conditions for the existence of extensions as in Theorem~\ref{thm:target_ext}. In view of the very recent results in \cite{VS24}, one could conjecture that the lower homotopy groups have to be merely finite and do not need to be trivial. 


\appendix

\addtocontents{toc}{\protect\setcounter{tocdepth}{1}}

\setcounter{equation}{0}
\renewcommand\theequation{A.\arabic{equation}}

\section{Retractions onto compact manifolds}\label{sec:retract_lem}

\subsection*{Proof of Lemma~\ref{lem:lip_retract}}
Recall that by assumption $N_\delta \subset Q_R$, where we defined $Q_R = (-R,R)^l$. Note that $N$ and \(\overline{N_\delta}\) share the same homotopy groups so that, in particular, \eqref{eq:pi} holds true with \(N\) replaced by \(\overline{N_\delta}\). Let now \(T\) be a triangulation of \(Q_R\setminus N_\delta\), and for \(i\in \{0,\dots, j\}\) denote by \(T^i\) its \(i\)-dimensional skeleton, that is, the union of all faces of the triangulation with dimension at most \(i\). 
The retraction \(P:Q_R\setminus X\to N\) is then obtained as the concatenation 
\begin{equation*}
P:\, Q_R\backslash X 
\overset{r^s}{\longrightarrow} T^{j+1}\cup \overline{N_\delta}
\overset{r}{\longrightarrow} \overline{N_\delta}
\overset{\xi}{\longrightarrow} N
\end{equation*}
with 
\begin{itemize}
	\item
	\(\xi\in C^1(\overline{N_\delta};N)\) the nearest point projection onto \(N\) from Lemma~\ref{lem:nearest_pp}, 
	\item
	\(r\in C^{0,1}(T^{j+1}\cup \overline{N_\delta};\overline{N_\delta})\) a Lipschitz retraction onto \(\overline{N_\delta}\),
	\item
	\(r^s\in C_{\loc}^{0,1}(Q_R\backslash X ;T^{j+1}\cup \overline{N_\delta})\) a locally Lipschitz retraction onto \(T^{j+1}\cup \overline{N_\delta}\).
\end{itemize}
In what follows we explain the construction of the maps \(r\) and \(r^s\). Results \ref{item:retraction} and \ref{item:inclusion} of Lemma~\ref{lem:lip_retract} are simple by-products of this construction. The claimed regularity of the retractions follows straightforwardly from their definition. Afterwards, we comment on the estimate in \ref{item:projection_grad_lp}.

\subsubsection*{The Lipschitz retraction \(r\)}
The construction follows an inductive argument. Choose some \(i\in \{0,\dots , j\}\) and some arbitrary simplex \(S\) in \(T^{i+1}\). For simplicity, since both shapes are homotopic, we assume that \(S=B^{i+1}\) is the unit ball of dimension \(i+1\) centered at the origin. The face of \(S\) is then given by its boundary \(\p S=\Sp^{i}\). Now, the task of constructing a retraction \(r: T^{j+1}\cup \overline{N_\delta}\to \overline{N_\delta}\) can be equivalently interpreted as the endeavor to extend the identity map on \(\overline{N_\delta}\) to \(T^{j+1}\cup \overline{N_\delta}\). Hence, arguing inductively, we suppose that the desired retraction is already given on the level of \(\p S=\Sp^{i}\), and we would like to extend the map to \(S=B^{i+1}\). This is possible because, as mentioned in Section~\ref{subsec:homotopy}, the pair \((B^{i+1},\sph^i)\) enjoys the {\it homotopy extension property}. This means that we can extend any given continuous map from \(\sph^i\) to \(B^{i+1}\) as soon as this holds true for a continuous function homotopic to this map on the level of \(\sph^i\). Now, firstly, any continuous map from \(\Sp^i\) to \(\overline{N_\delta}\) is contractible to a single point, thus homotopic to it, since \(\pi_i(N)=0\) by assumption \eqref{eq:pi}; and secondly, any constant map on \(\Sp^i\) can be extended constantly to the inside \(B^{i+1}\). This implies that, in fact, any map \(\p S\to \overline{N_\delta}\) can be extended to some map \(S\to \overline{N_\delta}\).

The above sketch can be made rigorous by using {\it Eilenberg's extension theorem}, which is a result in cohomology theory (cf.\ \cite[Chapter~6, Section~6]{Hu59}; see also the proof of Lemma~6.1 in \cite{HaLi87}). 

\subsubsection*{The locally Lipschitz retraction \(r^s\)}
So far, the retraction has been defined on the set \(T^{j+1}\cup \overline{N_\delta}\). We would like to proceed as outlined above, and keep lifting the retraction to higher and higher simplices by extending the given map on the faces to the inside of the simplices. Unfortunately, we have already exhausted the connectedness properties of $N$ (hence of $N_\delta$) up to order $j$, so for this procedure we can no longer rely on assumption \eqref{eq:pi}. Instead, we proceed in a more direct way. Assume that we are given some map \(f\) on some simplex boundary \(\p S=\Sp^{j+1}\) with values in \(T^{j+1}\cup \overline{N_\delta}\). We define the extension of such a map to \(S=B^{j+2}\) by
\begin{equation}
	f^s(x)
	\coloneqq f\left(\frac{x}{|x|}\right)
	\quad \forall x\in B^{j+2}.
\end{equation}
The price to pay is that \(f^s\) develops a point-singularity at the center of the ball. Repeating this construction for \(f\) the identity map and any simplex in \(T^{j+2}\), and then following the same steps to obtain a retraction from the next higher skeleton \(T^{j+3}\) onto \(T^{j+1}\cup \overline{N_\delta}\), results in a one dimensional singularity connecting the point-singularities emerging on the level of \(T^{j+2}\). Proceeding like this up to \(T^l=Q_R\) eventually leads to the retraction \(r^s:Q_R\backslash X \to T^{j+1}\cup \overline{N_\delta}\), where \(X\) is a (\(l-j-2\))-dimensional dual skeleton with respect to the triangulation \(T\) (see the proof of Lemma~2.2 in \cite{BPV14}). In particular, the singularity set \(X\) is contained in a finite union of (\(l-j-2\))-dimensional planes.

Concerning the estimate in Lemma~\ref{lem:lip_retract}\,\ref{item:projection_grad_lp}, first note that by \cite[Lemma~2.2]{BPV14} there holds
\begin{equation}\label{eq:basic_retract_est}
	|DP(y)|
	\le \frac{C}{\dist(y,X)},
\end{equation}
where the constant depends on \(l\) and \(N\). Hence, remembering that \(X\) is contained in a finite union of (\(l-j-2\))-dimensional planes, it will be sufficient to estimate the \(L^p\)-norm of the expression in \eqref{eq:basic_retract_est} where we assume that \(X\) is already given by one of these planes. Up to a unitary rotation, we might assume that we can decompose \(y\in Q_R\) as \(y=(y',y'')\) with \(y'\in X\) and \(y''\in X^{\perp}\), for \(X^\perp\) the orthogonal complement of \(X\) in \(Q_R\). We compute with the Fubini theorem
\begin{align}
\begin{split}
	\int_{Q_R} |DP(y)|^p {\dl y}
	&\le C \int_{Q_R} \dist(y,X)^{-p} {\dl y}\\
	&= C \int_X\int_{X^\perp} |y''|^{-p} {\dl y''}{\dl y'}
	\le C\mathcal{H}^{l-j-2}(X)\int_{X^\perp} |y''|^{-p} {\dl y''}.
\end{split}
\end{align}
The right-hand side, and in particular the integral, are now bounded by some constant depending on \(l\), \(N\), and \(R\), because \(X^{\perp}\) has dimension \(j+2\) and \(p<{j+2}\) by assumption. This concludes the proof of Lemma~\ref{lem:lip_retract}. 


\section*{Acknowledgements}
This work was primarily initiated and conducted while CG was affiliated with TU Wien, supported by the Austrian Science Fund (FWF) through the projects \href{https://doi.org/10.55776/F65}{10.55776/F65}, \href{https://doi.org/10.55776/V662}{10.55776/V662}, and \href{https://doi.org/10.55776/Y1292}{10.55776/Y1292}. For the final stages of this research, CG acknowledges support from the European Union's Horizon Europe research and innovation programme under the Marie Sk\l odowska-Curie grant agreement No 101102708.
The research of LH has been funded by the Austrian Science Fund (FWF) through the project \href{https://doi.org/10.55776/P34609}{10.55776/P34609}.
LH also acknowledges support from Grant PCI2024-155023-2 funded by\linebreak MCIN/AEI/10.13039/501100011033 and the European Union. 
The research of VP has been funded by the Austrian Science Fund (FWF) through the projects \href{https://doi.org/10.55776/I4052}{10.55776/I4052} and \href{https://doi.org/10.55776/Y1292}{10.55776/Y1292}.
CG and VP also acknowledge support from the Austrian Fe\-de\-ral Ministry of Education, Science and Research (BMBWF) through the OeAD-WTZ project CZ09/2023.
CG and VP are members of the GNAMPA group of INdAM (Istituto Nazionale di Alta Matematica "Francesco Se\-ve\-ri") and
acknowledge support from the project CUP\_E53C22001930001.


\bibliography{bibliographie}
\bibliographystyle{abbrvurl}

\end{document}